\documentclass[reqno,a4paper,11 pt]{article}
\usepackage{amssymb}
\usepackage{amsmath}
\usepackage[margin=0.7in]{geometry}
\usepackage{latexsym,amsmath}
\usepackage{amsmath}
\usepackage{amsfonts}
\usepackage[english]{babel}
\usepackage[dvips]{graphicx}
\usepackage{graphics}
\usepackage{makeidx}
\usepackage[dvips]{color}

\usepackage{amsmath,amssymb}
\def \E {\mathbb{E}}
\def \P {\mathcal{P}}
\def \e {\varepsilon}

\def \R {\mathbb{R}}
\newtheorem{theorem}{Theorem}[section]
\newtheorem{lemma}[theorem]{Lemma}

\newtheorem{hypothesis}[theorem]{Hypothesis}
\newcommand{\rmd}{\mathrm{d}}
\newtheorem{example}[theorem]{Example}

\newtheorem{remark}[theorem]{Remark}
\numberwithin{equation}{section}
\def\Dim{\noindent\hbox{{\bf Proof.}$\;\; $}}          %     end of proof
\def\finedim{{\hfill\hbox{\enspace${ \square}$}} \smallskip}    %     empty square
\def\sqr#1#2{{\vcenter{\vbox{\hrule height .#2pt
     \hbox{\vrule width .#2pt height#1pt \kern#1pt \vrule
     width .#2pt} \hrule height .#2pt}}}}
\def\square{\mathchoice\sqr54\sqr54\sqr{4.1}3\sqr{3.5}3}

\begin{document}
\thispagestyle{empty}
\parindent=0pt

\title{Singular Limit of BSDEs and Optimal Control \\ of two Scale Stochastic Systems in Infinite Dimensional Spaces}
\author{ %Francois Delarue\\ Laboratoire J.A.Dieudonn\'e \\
% UMR CNRS-UNS N°7351\\
% Universit\'e de Nice Sophia-Antipolis\\
% Parc Valrose\\
% France-06108 NICE Cedex 2\\
 %{\tt e-mail: delarue@unice.fr}
 %\\  \\ 
Giuseppina Guatteri, \\
Dipartimento di Matematica, \\
Politecnico di Milano, \\
Piazza Leonardo da Vinci 32, \\
20133 Milano,
Italia. \\
 {\tt e-mail: giuseppina.guatteri@polimi.it} \\ \\
Gianmario Tessitore,\\
Dipartimento di Matematica e Applicazioni \\
Universit\`{a} Milano-Bicocca, \\
 via R. Cozzi 53 - Edificio U5,\\ 20125 Milano, Italia.\\
{\tt e-mail: gianmario.tessitore@unimib.it
 }} 

\date{}
\maketitle
\begin{abstract}
In this paper we  study by probabilistic techniques the convergence of the value function for a two-scale,  infinite-dimensional, stochastic controlled system  as the ratio between the two evolution speeds diverges. 
The value function is represented as the solution of a \textit{backward stochastic differential equation} (BSDE) that it is shown to converge towards a \textit{reduced} BSDE. The noise is assumed to be additive  both in the slow and the fast equations for the state. Some non degeneracy condition on the slow equation is required. The limit BSDE involves the solution of an \textit{ergodic} BSDE and is itself interpreted as the value function of an auxiliary stochastic control problem on a reduced state space.
\end{abstract}

\section{Introduction} The purpose of this  paper is to give a representation of  the limit of the value functions of a sequence of optimal control problems for a singularly perturbed infinite dimensional state equation. Namely we consider the following system of controlled stochastic differential equations:
\begin{equation}\begin{cases}\label{controllointro}
dX^{\e,\alpha}_t= AX^{\e,\alpha}_t +  b(X^{\e,\alpha}_t,Q^{\e,\alpha}_t,\alpha_t)dt \,+ Rdl{W}^1_t, & X_0=x_0, \\ 
 \e dQ^{\e,\alpha}_t= (BQ^{\e,\alpha}_t+ F(X^{\e,\alpha}_t,Q^{\e,\alpha}_t) \,dt +  G\rho(\alpha_t)dt+ \e^{1/2} G \, d {W}^2_t, & Q_0=q_0,
\end{cases}
\end{equation}
where both state components $X^{\e,\alpha}$ and $Q^{\e,\alpha}$  take values in an Hilbert space. In the above equation $A$ and $B$ are unbounded linear operators, $\alpha$ represents the control, $  ({W}^1_t)_{t\geq 0}$,  $  ({W}^2_t)_{t\geq 0}$ are infinite dimensional cylindrical Wiener processes, $b$, $F$, $\rho$ are functions  and $R$ and $G$ are bounded linear operators satisfying suitable assumptions. We notice that the presence of the constant $\e$ in the second equation corresponds to the fact that $Q$ evolves with a speed which is larger by a factor $1/\e$ then the speed of evolution of the component $X$. In other words the above equation is a good model for a so called \textit{two scale system}. The optimal control problem is then completed by a standard cost functional of the form
:
\begin{equation} \label{costointro}J^{\e}(x_0,q_0,\alpha):=\E\left( \int_0^1 l(X^{\e,\alpha}_t, Q^{\e,\alpha}_t, \alpha_t)dt + h (X^{\e,\alpha}_1)\right),\end{equation}
and the value function if defined in the usual way:
\begin{equation}\label{intro-vf}
V^{\e}(x_0,q_0):=\inf_{\alpha} J^{\e}(x_0,q_0,\alpha),
\end{equation}
where the infimum is extended over a suitable class of progressively measurable control processes $(\alpha)$.  

Our purpose is to give a characterization of the limit of $V^{\e}(x_0,q_0)$ as $\e$ (that is the ratio between the speed of slow and the quick evolution) converges to $0$.

\medskip 

Several authors have studied the  convergence of singular stochastic control problems in finite dimensional spaces, see for instance \cite{AlvBar},  \cite{Bensoussan}, \cite{KabPer}, \cite{KabRungg}, \cite{Kus}. In particular \cite{AlvBar} has been an inspiration for the present work. In that paper the authors represent the value function of a singular stochastic control problem, in finite dimensions, by the solution, in viscosity sense, of an Hamilton Jacobi Bellman equation.  Then they show, by PDE methods  their convergence towards the solution, again in viscosity sense, of a \textit{reduced} parabolic PDE with smaller state space and a new nonlinearity usually called \textit{effective Hamiltonian}. Such analysis is performed in the case of periodic boundary conditions. Although  PDE techniques  perfectly fit the finite dimensional case allowing to cover 
general situations (including state equations with control dependent diffusions) they seem not to be adaptable to the infinite dimensional case, and consequently to the case of two scale stochastic control problems for stochastic PDEs. The reason essentially is the difficulty of handling, by  analytic tools and viscosity solutions, parabolic equations in infinite variables. Namelly comparison of viscosity solution (and consequently their uniqueness) always require, in  infinite dimensional frameworks,  additional artificial assumptions (see for instance the requirement of $B$-continuity and of trace class noise in \cite{FGS} and \cite{Swi}) that would not allow to cover our  case (see, as well, the discussion in the Introduction of \cite{FuTes}).

\medskip

In this paper  we choose a completely different approach based on Backward Stochastic Differential Equations, BSDEs in short, (see  \cite {ParPeng}, and  \cite{FuTes} as a reference, respectively, for the finite and infinite  dimensional case) that has already proved to be well adapted to infinite dimensional extensions. This choice eventually allows us to give a representation of  the limit of $V^{\e}(x_0,q_0)$ (see \eqref{intro-vf}) in a general Hilbertian framework that constitutes, at our best knowledge, the first result in this direction. Moreover our assumptions are general enough to cover a pretty  large class of two scale systems of controlled partial differential equations, possibly driven by cylindrical Wiener processes (see,  for instance, the system of controlled reaction diffusion equations driven by space-time white noise in Example \ref{esempio}). 
As a counterpart we notice that we consider state equation in which the control only affects the drift and in which the noise of the slow component is assumed to be non-degenerate.

\medskip

We try now to give a few more details on our method and results. To start with  we consider,
 for each $\e >0$, the following uncontrolled {\em forward-backward} system:
\begin{equation}\label{SystemEpsilonIntro}
\begin{cases}
dX_t= AX_t  \, dt + R  \, dW^1_t, & \\ 
 \e dQ^\e_t= (BQ^\e_t+ F(X_t,Q^\e_t) \,dt + \e^{1/2}  \,G \, dW^2_t, &  \\ 
- dY^{\e}_t  =  \psi(X_t,Q^\e_t,Z^\e_t, \Xi^\e_t / \sqrt{\e}) \,dt - Z^\e _t d\, W^1_t  - \Xi^\e_t dW_t^2, \\ 
X_0= x_0 \quad Q^\e_0=q_0, \quad {Y}^\e_1=h(X_1),
\end{cases} 
\end{equation}
where $\psi$ will eventually be the Hamiltonian corresponding of the stochastic control problem:  $$\psi(x,p,z,\xi)=\inf_{\alpha \in U} \{l(x,q,\alpha)+ z[R^{-1}b(x,q,\alpha)] + \xi \rho(\alpha) \}. $$
Then, once we have a solution $(X,Y^\e,Z^\e )$ to system \eqref{SystemEpsilonIntro}, we exploit the well known identification between $Y^\e_0$ and $V^{\e}(x_0,q_0)$
(see \cite{ElKPeQu} or \cite{FuTes})  in order to study the limit of the value functions by the limit of the sequence $Y^\e_0$ as $\e \rightarrow 0$.
 Our main result is indeed stated in terms of $Y^\e$, that is, see Theorem \ref{main}, we  prove that:
$$Y^\e_0\rightarrow \bar Y_0, \qquad \mathbb{P} -a.s.$$
where $(X,\bar Y,\bar Z)$ is the unique solution of the following decoupled forward backward system of stochastic differential equations:
$$
\begin{cases}
d X_t =   AX_t \, dt + R\, d W^1 _t, \\
-d\bar{Y}_t  = \lambda({X}_t,\bar{Z}_t)\,dt - \, \bar{Z}\, dW^1_t,   \\
X_0= x_0, \quad \bar{Y}_1=h(X_1).
\end{cases} 
$$
The statement of the above mentioned result is formulated and proved in Section 5 as a general result on singular limits of  BSDEs since it is independent of its control theoretic interpretation and, we believe, the proving argument has some interest on its own.
It is worth mentioning
 that the `reduced nonlinearity' $\lambda$ is itself a component of the unique solution $(\check{Y}, \check{Z}, \lambda)$ of the parametrized version of a, so called,  \textit{Ergodic} BSDE (see \eqref{BackInfHor} and Theorem \ref{Teolambda}) similar the ones  introduced in \cite{FuHuTes} (see \cite{DebHuTess} and \cite{HuRichou} as well). Function $\lambda$  can also be interpreted as  the optimal cost of an ergodic optimal control problem, see Remark \ref{rem-lambdacontrol}. Moreover, as it happens in  the finite dimensional case, the space in which the above reduced forward-backward system lives is a subspace of the original one (corresponding to the slow evolution alone).
 
 As a by-product of our main result, using the Bismut Elworthy formula in \cite{FuTes_BE} we immediately get that the solution of the reduced BSDE, and therefore the limit value function, depends on $x_0$ in a differentiable way  and is linked to the unique \textit{mild} solution of a semilinear parabolic PDE in infinite dimensional spaces:
$$ \left\{\begin{array}{l}\displaystyle
\frac{\partial v(t,x)}{\partial t}+\frac{1}{2}\hbox{Tr}[RR^*\nabla^2_x  v(t,x)] = \lambda (
x,\nabla v(t,x)),\qquad t\in [0,1],\,
x\in H,\\
\displaystyle v(1,x)=h(x).
\end{array}\right.
$$
Finally, in the last section, exploiting the concavity of $\lambda$ we give a representation  of $\bar{Y}_t$ as the value function of an auxiliary stochastic control problem on a reduced state space.

\smallskip
The paper is organized in the following way. In Section 2 we set the notation and we introduce some  functional spaces  while Section 3 contains some estimates on the two scale state equation that will be useful in the paper.
In Section 4 we introduce parametrized ergodic BSDEs and study their regularity with respect to parameters. In Section 5 we state the form of the limit equations and prove a convergence result for BSDEs that represents  the main technical issue of this paper. In Section 6, we finally link our results to the stochastic singular control problem. Finally, in section 7 we interpret the solution of the reduced BSDE in terms of a stochastic optimal control problem.

\section{Notation }
Given a Banach space $E$,
the norm of its elements $x$ will be denoted
by $|x |_E$, or even by $|x|$ when no confusion is possible. If $F$ is
another Banach space, $L(E,F)$ denotes the space of bounded linear
operators from $E$ to $F$, endowed with the usual operator norm. When $F=\mathbb{R}$ the dual space $L(E,\mathbb{R})$ will be denoted by $E^*$.  
The letters $\Xi$, $H$ and $K$ will always be used to denote Hilbert spaces.
The scalar product is denoted $\langle \cdot, \cdot \rangle$, equipped with a
subscript to specify the space, if necessary. All  Hilbert
spaces are assumed to be real and separable and  the dual of a Hilbert space will never be identified with the space itself. By $L_2(\Xi,H)$ and
$L_2(\Xi,K)$
 we denote the spaces of Hilbert-Schmidt operators from $\Xi $ to $H$
 and to $K$, respectively. Finally  $\mathcal{G}(K,H)$ is the space of all  Gateaux differentiable mappings $\phi$
 from $K$ to $H$ such that the map $(k,v)\rightarrow \nabla \phi(k)v$ is continuous from $K\times K$ to $H$; see \cite{FuTes} for details.
 \vspace{5pt}

Let $W^1 = (W^1_t)_{t\geq 0}$ and $W^2 = (W^2_t)_{t\geq 0} $ be two independent cylindrical Wiener processes with
values in $\Xi$, defined on a
complete
probability space $(\Omega, \mathcal{F}, \mathbb{P})$.
By $\{\mathcal{F}_t, \ t \in [0,T] \}$ we will denote the natural filtration
of $(W^1,W^2)$, augmented with the family $\mathcal{N}$ of
$\mathbb{P}$- null sets of $\mathcal{F}$. Obviously, the filtration
$(\mathcal{F}_t)$ satisfies the usual conditions of right-continuity and completeness. All the concepts
of measurability for stochastic processes will refer to this filtration.
By $\mathcal{P}$ we denote the predictable $\sigma$-algebra on
$\Omega \times [0,T]$ and by $\mathcal{B}(\Lambda)$ the Borel
$\sigma$-algebra of any topological space $\Lambda$.

Next we define  the following two classes of stochastic processes with values in
a Hilbert space $V$.  Given an arbitrary time horizon  $T$ and constant $p\geq 1$:
\begin{itemize}
\item $L^p_\P (\Omega\times [0,T];V)$ denotes the space of
equivalence classes of processes $Y \in L^p (\Omega\times
[0,T];V)$ admitting a predictable version. It is endowed with the norm
\[ |Y|= \Big(\E \int_0^T |Y_s|^p \, ds\Big)^{1/p}. \]
\item $L^{p, loc}_\P (\Omega\times [0,+\infty[;V)$ denotes the set of processes defined on $\mathbb{R}^+$ such that their  restriction to an arbitrary  $[0,T]$ belongs to $L^p_\P (\Omega\times [0,T];V) $.  
   \item $L^p_{\mathcal{P}}(\Omega;C([0,T];V))$
     denotes the space of
    predictable processes $Y$ with continuous paths in $V$, such
    that the norm
    \[  \|Y\|_p  = (\E \sup _{s \in [0,T]} |Y_s|^p)^{1/p}\]
    is finite. The elements of $L^p_{\mathcal{P}}(\Omega;C([0,T];V))$
    are identified up to indistinguishability.
    \item $L^{p, loc}_\P (\Omega;{C} [0,+\infty[;V))$ denotes the set of processes defined on $\mathbb{R}^+$ such that their  restriction to an arbitrary  $[0,T]$ belongs to $L^p_{\mathcal{P}}(\Omega;C([0,T];V))$. 
\end{itemize}
Given $\Phi$ in  $L^2_\P (\Omega\times
[0,T];L_2(\Xi,V))$, the It\^o stochastic integrals $\int_0^t
\Phi_s \,dW^1_s$ and  $\int_0^t \Phi_s \,dW^2_s$, $t \in
[0,T]$, are  $V$-valued martingales belonging to
$L^2_{\mathcal{P}}(\Omega;C([0,T];V))$.
% The previous definitions
%have obvious extensions to processes defined on  the entire positive real line $\R^+$.

%%%%%%%%%%%%%%%%%%%%%%%%%%%%%%%%%%%%%%%%%%%%%%%%%%%%%%%%%%%%%%%%%
\section{The forward system}
%%%%%%%%%%%%%%%%%%%%%%%%%%%%%%%%%%%%%%%%%%%%%%%%%%%%%%%%%%%%%%%%%
For arbitrarily fixed $x_0 \in H$ and $q_0 \in K$ we consider
 the following  system of  stochastic differential
equations in $H \times K$:
\begin{equation}
\label{eqX}
\begin{cases}
dX_t= AX_t dt +  R\, dW^1_t, & X_0=x_0,\; t\geq 0, \\ \\
 \e dQ^\e_t= (BQ^\e_t+ F(X_t,Q^\e_t) )\,dt + \e^{1/2} G \, dW^2_t, & Q^\e_0=q_0,\; t\geq 0,
\end{cases}
\end{equation}
 where the ``slow'' variable $X$ takes its values
in $H$ and the ``fast'' variable $Q^\e$ takes its values in $K$,
$\e \in ]0,1]$ is a small parameter. 

Finally $A: D(A)\subset H \to
H$ and $B: D(B) \subset K \to K$ are unbounded linear operators
generating $C_0$- semigroups $\{ e^{tA} \}_{t \ge 0}$ and $\{
e^{tB} \}_{t \ge 0}$ over $H$
and $K$, respectively, while  $R$ and $G$ are linear bounded operators from $\Xi$ to $H$ (respectively to $K$). \\
Moreover, we make the following, standard assumptions:
\begin{hypothesis}\label{A.1}
$A: D(A)\subset H \to H$
is a 
linear, unbounded  operator that generates a  $C_0$- semigroup $\{
e^{tA} \}_{t \geq 0}$, such that  $|e^{tA}|_{L(H,H)} \le M e ^{ \omega_A t},  t \geq 0$ for some positive constants $M_A$ and $\omega_A$.
$B: D(B)\subset K \to K$ is a 
linear, unbounded  operator that generates a  $C_0$- semigroup $\{ e^{tB} \}_{t \ge 0}$ such that $|e^{tB}|_{L(K,K)} \le   M_B e ^{ \omega_B t}, t \geq 0$
 for some $ M_B, \omega_B>0$. \\
Moreover there exist  constants $L>0$ and $\gamma \in [0,\frac{1}{2}[$ s.t.:
\begin{align}
&|e^{sA}|_{L_2(\Xi,H)}+ |e^{sB}|_{L_2(\Xi,K)}\nonumber \leq L
s^{-\gamma},\quad \forall s\in [0,1].\nonumber 
\end{align}

\end{hypothesis}
\begin{hypothesis}\label{A.2}
$F:  H
\times K \to K$ is bounded and there exists a constant $L_F$ for which:
\begin{equation}
| F(x,y)-
F(u,v)| _K \leq L_F( |x-u|_H + |y-v|_K).\nonumber 
\end{equation}
 for every $x,u
\in H$, $y,v \in K$.

Moreover we assume that for every $x \in H$, $F (x,\cdot)$ is Gateaux differentiable, more precisely, $F (x,\cdot) \in \mathcal{G}^1(K,K)$.

\end{hypothesis}
\begin{hypothesis}\label{A.3} $B +F$ is  dissipative i.e.
there exists some $\mu >0$ such that:
\[ \langle B q+F (x,q) - (Bq' +F(x,q')), q-q' \rangle
 \leq -\mu |q-q'|^2, \]
for all $
x\in H,
 q,q' \in D(B)$.
\end{hypothesis}
\begin{hypothesis}
 \label{A.4}  $R \in L(\Xi;H)$, $G \in {L}(\Xi;K)$ and  moreover $R$ admits a bounded right inverse $R^{-1}\in L(H;\Xi)$.
\end{hypothesis}
%\begin{remark}\label{dissip}
%From the assumptions \eqref{A.1}, \eqref{A.2} and %\eqref{A.3} it follows that $B-L_F I$ is strongly %dissipative with
%dissipativity constant $\mu$.
%Indeed for every $y \in D(B)$:
%\begin{align*}
%&\langle (B-L_F I) y,y\rangle =\langle (B+ F(x^0,%\cdot)) y -  (B+ F(x^0,\cdot))0,y \rangle
%\\ & + \langle (F(x^0,0) - F(x^0,y)),y \rangle - %L_F |y|^2  \leq -\mu |y|^2
% \end{align*}
%Moreover, \eqref{A.2} and \eqref{A.3} imply for $B$ assumption \eqref{A.1} with $M_B=1$ and $\omega_B= L_F - \mu.$
%Similarly we have that $F+ L_F$ is dissipative. 
%\end{remark}
Given any cylindrical Wiener process $(\beta_t)_{t\geq 0}$ with values in $\Xi$ we denote by $(\beta^B_t)_{t\geq 0}$ the stochastic convolution
$$\beta^B_s=\int_0^s e^{(s-\ell)B}Gd\beta_{\ell}.$$
In the following we shall assume, as in \cite{FuHuTes}, that:
\begin{hypothesis}
 \label{A.5}  $\sup_{s>0}\E|\beta^B_s|^2<\infty$.

\end{hypothesis}

\begin{remark}
Notice that since $(\beta_t)$ is a centered gaussian process this implies that,   $\forall p\geq 1$ it holds $\sup_{s>0}\E|\beta^B_s|^p<\infty$. Moreover hypothesis $\ref{A.5}$ is  verified whenever $B$
is a strongly dissipative operator.
\end{remark}
We collect here two results we will use in the sequel.
We do not provide the proof of the first, that can be found for instance in \cite[Proposition 3.2]{FuTes}. Regarding the second result, for the reader's convenience, we briefly
 report the argument which is a slight modification of the one in \cite[section 6.3.2.]{daza2}.

\begin{lemma}\label{stimaX} Under Hypothesis \ref{A.1}  and \ref{A.4} the slow equation in system  (\ref{eqX}) admits a unique mild solution $X^{x_0}_t$ that has continuous trajectories and  
 for all $ p \geq 1$ satisfies:
 \begin{equation}\label{stimaslow}
 \E(\sup_{t\in [0,1]}|X^{x_0}_t|^p)\leq c_p(1+|x_0|^p), \qquad x_0 \in H,
 \end{equation}
 for some positive constant $c_p$ depending only on $p$ and on the quantities introduced in the hypotheses.
\end{lemma}

\begin{lemma}\label{dipendcont}

Let $(\Gamma_s)_{s\geq 0}$ be a  given, $H$-valued, predictable process with $\Gamma \in L^{p,loc}_{\mathcal{P}}(\Omega\times[0,\infty[;H)$ and let $(g)_{s\geq 0}$ be a given, $K$-valued, process with $g \in L^{p,loc}_{\mathcal{P}}(\Omega\times[0,+\infty[;K)$ for some $p\geq 1$.

Then the  following  equation:
\begin{equation} \label{problemalimiteyX}
d{Q}_s=(B{Q}_s+ F(\Gamma_s,{Q}_s))\,ds +g_sds+ Gd\beta_s, \ s \geq 0,
\hspace{20pt} {Q}_0= q_0,
\end{equation}
 admits a unique mild solution $Q \in L^{p,loc}_{\mathcal{P}}(\Omega;C([0,+\infty[;K))$.

Under hypotheses \eqref{A.1}--\eqref{A.5}, there exists a constant $k_p$ (independent on $T$) such that for all $T>0$:
\begin{equation}
\sup_{s\in [0,T]}\E |{Q}_s|^p \leq k_p(1+ |q_0|^p + \sup_{s\in [0,T]}\E|\Gamma_s|^p+\sup_{s\in [0,T]}\E|\beta^B_s|^p+\sup_{s\in [0,T]}\E|g_s|^p).
\end{equation}
Moreover if $(\Gamma'_s)_{s\geq 0}$ is another  $H$-valued, predictable processes in $ L^{p,loc}_{\mathcal{P}}(\Omega\times[0,\infty[;H)$ and $Q'$ is the mild solution of equation:
 $$
d{Q'}_s=(B{Q'}_s+ F(\Gamma'_s,{Q'}_s))\,ds +g_sds+ Gd\beta_s, \ s \geq 0,
\hspace{20pt} {Q}_0= q_0,
$$
then, for all $T>0$,
 $$|Q_T-Q'_T|\leq K \int_0^T e^{-\mu(T-\ell)}|\Gamma_{\ell}-\Gamma'_{\ell}| \, d\ell,\;\; \hbox{ $\mathbb{P}$-a.s., }$$
  where again $K$ does not depend on $T$.
\end{lemma}
\Dim 
Let $Z_s=e^{\mu s}(Q_s-\beta^B_s)$.  By Ito rule (going through Yosida approximations) we deduce that $Z$ is the mild solution of the following equation 
$$d {Z_s}=\mu Z_s\, ds+BZ_s\, ds+e^{\mu s}F(\Gamma_s,e^{-\mu s}Z_s+\beta^B_s))ds+e^{\mu s}g_sds.$$
Differentiating $\sqrt{|Z_s|^2+\e}$ (going, once more, through Yosida approximations), using dissipativity of $B+F$, see  hypothesis \ref{A.3}, we obtain
$$|Z_s|\leq \sqrt{|Z_s|^2+\e}\leq \sqrt{q^2_0+ \e}+\int_0^s e^{\mu \ell}\left|F(\Gamma_{\ell},\beta^B_{\ell})+g_{\ell}\right| d\ell +\mu \int_0^s \left[\sqrt{|Z_{\ell}|^2+\e}-|Z_\ell| \right]d\ell.$$
Letting $\e \rightarrow 0$, by dominated convergence we obtain:
$$|Z_s|\leq |q_0|+\int_0^s e^{\mu \ell}\left|F(\Gamma_{\ell},\beta^B_{\ell})+g_{\ell}\right| d\ell. $$ 
Recalling the definition of $Z$
we conclude:
$$|Q_s|\leq |\beta^B_s|+ e^{-\mu t}|q_0|+\int_0^s e^{-\mu(s- \ell)}\left|F(\Gamma_{\ell},\beta^B_{\ell})+g_{\ell}\right| d\ell. $$
and by Holder inequality (for the last term):
$$|Q_s|^p\leq 3^p |\beta^B_s|^p+ 3^p e^{-p\mu t}|q_0|^p+ 3^p\left(\int_0^s e^{-p^*\frac{\mu}{2}(s- \ell)}d \ell\right)^{p/p^*} \int_0^s e^{-p\frac{\mu}{2}(s- \ell)}|F(\Gamma_{\ell},\beta^B_{\ell})+g_{\ell}|^p d\ell. $$
The claim then follows from Hypothesis \ref{A.2}.
%Taking into account \eqref{A.3} we have that, as %in  \cite{FuHuTes}, Proposition 3.6
%%\sup_{t\in [0,\infty)}\E |\hat{Q}^{\hat{X},q_0}%_t|^2 \leq K(1+ |q_0|^2 + \sup_{t\in [0,\infty)}\E|\hat{X}_t| ^2)
%\end{align*}
%where the constant $K$, thanks to the dissipativity condition, can be chosen independent of $T$.

%Using again the \eqref{A.3} in order to evaluate %the differential we end up with:
%\begin{align*}
%d_t |\hat{Q}^{\hat{X}^1,q_0}_t- \hat{Q}^{\hat{X}^2,q_0}_t|^2\leq  -\mu |\hat{Q}^{\hat{X}^1,q_0}_t- \hat{Q}^{\hat{X}^2,q_0}_t|^2\, dt + L_F  |\hat{X}^1_t- \hat{X}^2_t|^2\, dt
%\end{align*}
%Therefore,
%%|Q^{\hat{X}^1,q_0}_t- Q^{\hat{X}^2,q_0}_t|^2\leq L_F\int_0^t e^{-\mu (t-s)}  |\hat{X}^1_s- \hat{X}^2_s|^2\, ds
%\end{align*}
% and this yields to conclude:
%\begin{align*} \sup_{r \geq 0} |Q^{\hat{X}^1,q_0}%%\, \E \sup_{r \geq 0} |\hat{X}^1_r-\hat{X}^2_r|^2. \quad \mathbb{P}- a.s.
%\end{align*}
The proof of the last statement is similar (and easier) noticing that:
$$d_s(Q_s - Q'_s)=B(Q_s - Q'_s)ds+[F(\Gamma_s,Q_s)-F(\Gamma'_s,Q'_s)]ds, $$
and then arguing as before.
\finedim

\noindent
If we fix $x\in H$, $q_0\in K$, choose $g\equiv 0$  and make a change of time $ s \to \e s$, then the fast equation in system \eqref{eqX} becomes 
 \begin{equation} \label{problemalimitey}
d\hat{Q}_s=(B\hat{Q}_s+ F(x,\hat{Q}_s))\,ds + Gd\hat{W}^2_s, \ s\geq 0,
\hspace{20pt} \hat{Q}_0= q_0.
\end{equation}
where  $\hat{W}_s^{2}={\e}^{-1/2}{W}_{\e s }^{2}$
is a cylindrical Wiener process.
 So \eqref{problemalimitey} is a special case of \eqref{problemalimiteyX}, and Lemma \ref{dipendcont}
applies. 

 We will denote by $\hat{Q}^{x,q_0}_s$ the unique mild solution of equation \eqref{problemalimitey}.

%\textbf{[Questo serve?]}
%\begin{lemma}\label{fastrescaled}
%For every fixed $x \in H$  consider the problem,
%then under hypotheses \ref{A.1}-\ref{A.2}-\ref{A.3}, \begin{equation} \label{problemalimiteybis}
%d\hat{Q}_s=(B\hat{Q}_s+ F(x,\hat{Q}_s))\,ds + d\hat{W}^2_s, \ s\geq 0,
%\hspace{20pt} {Q}_0= q_0.
%\end{equation}
%Then for  every $q$ there exists a unique mild solution $\hat{Q}^{x,q_0}$ and  there is a unique invariant measure $\mu^x(dy)$ for which the following convergence result holds true for every $\phi: K\to \R$, $\phi$ Lipschitz:
%\begin{displaymath}\label{convesponen}
%|P_t\phi(x,q_0)- \langle \phi,\mu^x\rangle| \leq C(1+ |x|+ |q_0|)\exp^{-\mu t /2}\| \phi \|_{Lip}, \ \qquad 
%\end{displaymath}
%where $P_t\phi(x,q_0)= \E(\phi(\hat{Q}^{x,q_0}_t))$ and $ \langle \phi,\mu^x\rangle= \int_K \phi(x,q_0)\, d\mu^x(dq_0)$ and $C$ is a constant independent on $T$.
%Moreover, for some positive constant $C$ depending on $\mu, L_F$, we have that
%\begin{equation}\label{LipMisura}
%| \langle \phi(\cdot),\mu^x \rangle- \langle \phi(\cdot),\mu^y \rangle | \leq  C || \phi ||_{Lip} |x-y|
%\end{equation}
%\end{lemma}
%\Dim
%The proof can be found in \cite{daza2}, Theorem 6.3.3, the dependence on the parameter $x$ is straightforward having all the coefficients sublinear growth with respect to $x$.
%\cite{DebHuTess}, Theorem 2.4 and the Appendix A.2..  
%\finedim
%%%%%%%%%%%%%%%%%%%%%%%%%%%%%%%%%%%%%%%%%%%%%%%%%%%%%%

\section{The ergodic BDSE parametrized}
We introduce a function $\psi: H\times K\times \Xi^* \times \Xi^* \to \R$. We will eventually (see Section 6) choose as $\psi$ the Hamiltonian of our control problem. Here we only assume that $\psi$ satisfies the following:
\begin{hypothesis}\label{B.3} Function $\psi$ is   measurable and there exist $L_q,  L_x, L_z ,L_\xi >0$ such that $\forall  \, q,q' \in K$,  $x, x' \in H$, $ \xi, \xi', z, z' \in \Xi^*$:
 \begin{multline*}|\psi(x, q,z ,\xi) - \psi(x', q', z', \xi')| \leq L_x (1+|z|)|x-x'| +  L_z |z-z'|+ L_q(1+|z|)|q-q'| +  L_\xi |\xi-\xi'|.%\\
%| \psi(x, q,z ,\xi)| \leq L(1+|z|+|q|)\\ 
\end{multline*}
Moreover we assume that  $\sup_{x\in H,q\in K} |\psi (x,q,0, 0)|<+\infty $
\end{hypothesis}

%Let us fix $x \in H$, $s\geq 0$ and ${q}_0 \in K$ 
%and consider of \eqref{problemalimitey}  starting at $s$ in $q$:
%\begin{equation} 
%d\hat{Q}_t=(B\hat{Q}_t+ F(x,\hat{Q}_t))\,dt +  G %d\hat{W}^2_t, \ t \geq s,
%\hspace{20pt} {Q}_s={q_0},
%\end{equation}
%and we denote by $\hat{Q}^{x,q_0}$ the solution % of \eqref{problemalimiteyX}.
The next result states existence of a solution to the so called {\em ergodic backward stochastic differential equation} (EBSDE):
\begin{equation}
\label{BackInfHor}
- \rmd \check{Y}_t=[\psi(x,\hat{Q}^{x,q_0},z,\check{\Xi}_t)-\lambda(x,z)]\,\rmd t - \check{\Xi}_t \rmd \hat{ W}_t^{2}, \quad \forall \, t\geq s.
\end{equation}
\begin{theorem}\label{Teolambda}
Under hypotheses  \ref{A.1},  \ref{A.2},  \ref{A.3}, \ref{A.4}, \ref{A.5} and  \ref{B.3} there exist measurable functions $\check{v}: H\times K \times \Xi^* \rightarrow \mathbb{R}$,
$\check{\zeta}: H\times K \times \Xi^*\rightarrow \mathbb{R}$, $\lambda: H \times \Xi^*\rightarrow \mathbb{R}$  with:
\begin{equation}\label{stimacheckv}|\check{v}(x,q,z)| \leq c(1+|z|)|q|,\end{equation}
(where  $c>0$  depends only on the constants introduced in the above mentioned hypotheses) such that the following holds: if we set: 
\begin{equation}\label{defdicheckY} \check{Y}^{x,q_0,z}_t= \check{v}(x,\hat{Q}^{x,q_0}_t,z), \quad \check{\Xi}^{x,q_0,z}_t= \check{\zeta}(x,\hat{Q}^{x,q_0}_t,z),\end{equation}
then  $\check{\Xi}^{x,q_0,z}$ is in $L^{2,loc}_\P([0,+\infty [,\Xi^*)$ and  $\mathbb{P}$-a.s. the EBSDE \eqref{BackInfHor} is safisfied by  $(\check{Y}^{x,q,z}_t, \check{\Xi}^{x,q,z}, \lambda(x,z))$ for all $ 0\leq t \leq T$.

%Thus, see for instance \cite{HuRichou}, 
%\begin{equation}
%\sup _{ t\geq 0}\E |\hat{Y}^{x,q,p}_t| \leq  %C (1+ |q|)
%\end{equation}
%for some constant  $C$ depending only on $ %M_B, \gamma $ and $\sup_{q \in K} |F(x,q)|$.

Moreover we have:
\begin{equation}\label{LipLambda}
|\lambda(x,z) -\lambda( x',z')| \leq {L^1_x} (1+|z|)|x-x'| +{L^1_z}|z-z'|,
\end{equation}
for some posive constants $L^1_x$ and $L^1_z$.
\end{theorem}
\Dim 
Fix $x\in H$ and $z\in \Xi^*$. In 
 \cite[Theorem 4.4  and Corollary 5.9]{FuHuTes} authors prove existence of functions  $\check{v}(x,\cdot\,,z)$, $\check{\zeta}(x, \cdot\,,z)$ and $\lambda(x,z)$ such that (\ref{stimacheckv}) holds and if  $\check{Y}^{x,q_0,z}$, $\check{\Xi}^{x,q_0,z}$ are defined as in (\ref{defdicheckY}), then $\check{\Xi}^{x,q_0,z}$ is in $L^{2,loc}_\P([0,+\infty[,\Xi^*)$ and $(\check{Y}^{x,q_0,z}_t, \check{\Xi}^{x,q_0,z}, \lambda(x,z))$  is a solution to equation   \eqref{BackInfHor}.

 Measurably of $\check{v}$, $\check{\zeta}$ and $\lambda$ with respect to all parameters follows by their construction (see again \cite{FuHuTes} Theorem 4.4).

$ $

We only need to prove \eqref{LipLambda}. 
Fixed $x,x' \in H$ and $z, z' \in \Xi^*$ we set $\widetilde{\lambda}=\lambda(x,z)-\lambda(x',z')$, $\widetilde{Y}= Y^{x,0,z} -Y^{x',0,z'},$ $ \widetilde{\Xi} =\check{\Xi}^{x,0,z} -\check{\Xi}^{x',0,z'},$ 
\begin{equation*}
\theta_t = \left\{ \begin{array} {ll}  \displaystyle \frac{\psi(x,\hat{Q}^{x,0}_r,z,\check{\Xi}^{x,0,z}_r)-\psi(x,\hat{Q}^{x,0}_r,z,\check{\Xi}^{x',0,z'}_r) }{|\check{\Xi}^{x,0,z}_r-\check{\Xi}^{x',0,z'}_r|^2_{\Xi ^*}}(\check{\Xi}^{x,0,z}_r-\check{\Xi}^{x',0,z'}_r), &  \hbox{ if }\; \check{\Xi}^{x,0,z}_r\not=\check{\Xi}^{x',0,z'}_r  \\
0 & \text{elsewhere}
\end{array}
\right. 
\end{equation*}
and
\begin{align*}
f_t &= \psi(x,\hat{Q}^{x,0}_r,z,\check{\Xi}^{x',0,z'}_r)- \psi(x',\hat{Q}^{x',0}_r,z',\check{\Xi}^{x',0,z'}_r ).\end{align*}
Then we have
\begin{equation*}
\widetilde{Y}_0+\lambda T= \widetilde{Y}_T + \int_t^T  f_r \,\rmd r - \int_t^T\widetilde{\Xi}_r (\theta_t dt +\rmd \hat{W}_r^{2}), \quad \forall \,  T \geq  t\geq 0.
\end{equation*}
So, by Girsanov theorem (notice that $(\theta_t)$ is uniformly bounded), there exists a probability  $\widetilde{\mathbb{P}}$  (mean value denoted by  $\widetilde{\mathbb{E}}$) 
such that $\widetilde{W}_t=\int_0^t \theta_{\ell} d\ell +\hat{W}_t^{2}$, $t\geq 0$,  is a cylindrical Wiener process. Consequently:
\begin{equation*}
{\lambda}T= \widetilde{Y}_T- \widetilde{Y}_0  + \int_0^T  f_r \,\rmd r - \int_0^T\widetilde{\Xi}_r \rmd \widetilde{W}_r, \quad \forall \,  T \geq  t\geq 0
\end{equation*}
and consequently:
\begin{align}\label{stimaquasifinalelambda}
| \lambda|\leq T^{-1}|\widetilde{Y}_0| + T^{-1}\widetilde{\E}|\widetilde{Y}_T|+ T^{-1} \int_0^T  \widetilde{\E} |f_s| \, \rmd s .
\end{align}
Thanks to hypothesis \ref{B.3} we get that for all $t \geq 0$:
\begin{equation*}
|f_t| \leq L_x(1+|z|)|x-x'| + L_z |z-z'| + L_{q}(1+|z|)|\hat{Q}^{x,0}_t-\hat{Q}^{x',0}_t| , \qquad \mathbb{P}-a.s.
\end{equation*}
We notice that with respect to $(\widetilde{W}_t)$
processes $\hat{Q}^{x,0}$ and  $\hat{Q}^{x',0}$ satisfy respectively
 \begin{eqnarray*}
 d\hat{Q}^{x,0}_s=(B\hat{Q}^{x,0}_s+ F(x,\hat{Q}_s^{x,0}))\,ds + \theta_s ds+ Gd\widetilde{W}_s, \ s\geq 0,\\
  d\hat{Q}^{x',0}_s=(B\hat{Q}^{x',0}_s+ F(x',\hat{Q}^{x',0}_s))\,ds + \theta_s ds+ Gd\widetilde{W}_s, \ s\geq 0,
 \end{eqnarray*}
 and  Lemma \ref{dipendcont} yields
 $|\hat{Q}^{x,0}_s-\hat{Q}^{x',0}_s|\leq (K/\mu) |x-x'|$ thus:
 \begin{equation}\label{stimadif}
| f_t|\leq (L_x+L_q K/\mu)(1+|z|)|x-x'|+L_z|z-z'|,\qquad \mathbb{P}-a.s. \text{ for all }t \geq 0.
 \end{equation}
 From Lemma \ref{dipendcont}  we also have that for every $T \geq 0$  and every $p\geq 1$,
 \begin{equation}
\sup_{s\in [0,T]}\widetilde{\E} |\hat{Q}^{x,0}_s|^p \leq k_p(1+ |x|^p +\sup_{s\in [0,T]}\widetilde{\E}|\theta_s|^p), 
\end{equation}
and 
 \begin{equation}\label{stimadip}
\sup_{s\in [0,T]}\widetilde{\E} |\hat{Q}^{x',0}_s|^p \leq k_p(1+ |x|^p +\sup_{s\in [0,T]}\widetilde{\E}|\theta_s|^p).
\end{equation}
Since $\theta$ is uniformly bounded it holds: 
 $$\displaystyle \sup_{t\in [0,\infty[} \widetilde{\E}(|\hat{Q}^{x,0}_t|^p+|\hat{Q}^{x,0}_t|^p)<\infty, $$
thus, by
 \eqref{stimacheckv}, we get that:
  $$\displaystyle \sup_{t\in [0,\infty[} \widetilde{\E}(|\widetilde{Y}_t|)<\infty.$$ 
  
  Consequently $T^{-1}\widetilde{\E}(|\widetilde{Y}_T|)
 \rightarrow 0$ as $T\rightarrow \infty$ and the claim follows by (\ref{stimaquasifinalelambda}) and (\ref{stimadif}) letting $T\rightarrow \infty$.

\finedim

%\begin{remark}
%Notice that, as already pointed out in \cite{FuHuTes}, there's no need to require the non degeneracy  of the operator $G$ in front of the noise in equation \eqref{problemalimitey}.
%This is a feature of the probabilistic approach, that seems able to overcome the non degeneracy condition required in the approach with viscosity solutions, see \cite{AlvBar}.
%\end{remark}
\begin{remark}
If, fixed $x$ and $z$, one restricts the class of triples $(Y,\Xi,\lambda)$ where to find a solution to equation \eqref{BackInfHor}, asking that there must be a constant $c>0$ (that may depend on $q_0$, $x$ and $z$) such that  $|Y_t| \leq c (1+ |Q_t|) $   $\mathbb{P}$-a.s. for every $t \geq 0$ then, see \cite[Theorem 4.6]{FuHuTes}, 
 the third component $\lambda$  of the solution id uniquely determined.
\end{remark}
\section{Limit equation and convergence of singular BSDEs}
We've eventually  got to the  {\em forward-backward system} for  $t \in [0,1]$
\begin{equation}\label{SystemEpsilon}
\begin{cases}
dX_t= AX_t  \, dt + R  \, dW^1_t, & \\ 
 \e dQ^\e_t= (BQ^\e_t+ F(X_t,Q^\e_t) \,dt + \e^{1/2}  \,G \, dW^2_t, &  \\ 
- dY^{\e}_t  =  \psi(X_t,Q^\e_t,Z^\e_t, \Xi^\e_t / \sqrt{\e}) \,dt - Z^\e _t d\, W^1_t  - \Xi^\e_t dW_t^2, \\ 
X_0= x_0 \quad Q^\e_0=q_0, \quad {Y}^\e_1=h(X_1),
\end{cases} 
\end{equation}
that, as  we will see in the sequel, is also  associated to a controlled multiscale dynamics. 
Function  $h: H\rightarrow \mathbb{R}$ satisfies:
\begin{hypothesis}\label{B.4} $h$ is  Lipschitz continuous with constant $L >0$.
\end{hypothesis}

We have that:
\begin{theorem}\label{esistenceepsilon}
Assume \ref{A.1}--\ref{A.5}, \ref{B.3} and \ref{B.4}. For every $\e >0$  there exists a unique  5-tuple of processes $(X, Q^\e,Y^\e, Z^\e,\Xi^\e)$,  with $ X\!\in L^2_{\mathcal{P}}(\Omega;C([0,1];H))$, $Q^\e \!\in L^2_{\mathcal{P}}(\Omega;C([0,1];K))$, $Y^\e\! \in L^2_{\mathcal{P}}(\Omega;C([0,1];\R))$, $Z^\e \!\in L^2_\P (\Omega\times [0,1];\Xi^*)$ and  $\Xi^\e \!\in L^2_\P (\Omega\times [0,1];\Xi^*)$ such that $\mathbb{P}-a.s.$  the  system  \eqref{SystemEpsilon} is satisfied for all $ t \in [0,1]$.
\end{theorem}
\Dim The proof is contained in \cite[Propositions 3.2 and  5.2]{FuTes}, we just notice that
the system is decoupled, so once the forward equation is solved then it becomes a known process in the backward equation. 
\finedim

\medskip

The purpose of our work is to study the limit behaviour of $Y^\e$ as $\e $ tends to $0$.

We introduce the candidate limit equation, that turns out to be a forward-backward system on the {\em finite horizon} $[0,1]$ and on the reduced state space $H$.
\begin{equation}\label{LimitEquation}
\begin{cases}
d X_t =   AX_t \, dt + R\, d W^1 _t, \qquad t \in [0,1], \\
-d\bar{Y}_t  = \lambda({X}_t,\bar{Z}_t)\,dt - \, \bar{Z}_t\, dW^1_t,   \\
X_0= x_0, \quad \bar{Y}_1=h(X_1).
\end{cases} 
\end{equation}
where $\lambda$ is defined in Theorem \ref{Teolambda}.

 One has that
\begin{theorem} \label{esistencebar}Under Hypothesis \ref{A.1}---\ref{A.5}, \ref{B.3}  and \ref{B.4},
there exists a unique triplet of processes 
$(X,\bar{Y},\bar{Z})$ with $X\!\in L^p_{\mathcal{P}}(\Omega;C([0,1];H))$, $\bar{Y}\,\in L^p_{\mathcal{P}}(\Omega;C([0,1];\R))$, $\bar{Z}\!\in L^p_\P (\Omega\times [0,1];\Xi^*)$ that   fullfils system  \eqref{LimitEquation}, $\mathbb{P}-a.s.$  for every $ t\in [0,1]$.
\end{theorem}
\Dim Thank to the regularity of $\lambda$, see  (\ref{LipLambda}), the proof of existence and uniqueness of the solution to equation \eqref{LimitEquation} is standard (see, for instance \cite[Proposition 4.3]{FuTes}).
\finedim

We can now state our main result:
\begin{theorem}\label{main}
Under Hypothesis \ref{A.1}---\ref{A.5}, \ref{B.3}  and \ref{B.4}, the following holds for  $\bar{Y}$ and   $ Y^\e $ found in Theorem \ref{esistenceepsilon} and Theorem \ref{esistencebar}  respectively: 
\begin{equation}\label{Conv}
\lim_{\e \to 0} Y^\e _0 =\bar Y_0.
\end{equation}
\end{theorem}
\Dim
%Let us consider the first term of the sum at R.H.S. Notice that $X^\e$ and $\bar{X}$ coincide, being the solution of the same equation [QUI BISOGNA CONTROLLARE]. 
%as it is, bearing in mind that, by \eqref{Lippsi}, the following holds:
%\begin{equation} \label{limphi1}
%|\psi(\bar{X}_t,Q^\e_t,{Z}^\e_t,\Xi^\e_t/\sqrt{\e})-\psi(\bar{X}_t,Q^\e_t,Z^\e_t,\Xi^\e_t/\sqrt{\e})|\leq L |%{Z}^\e_t-\bar{Z}_t|
%\end{equation}
We start by noticing  that if we slow down time, that is, for $s\in [0,1/\e[$ we set:
 $\hat{Q}^{\e}_s=Q^{\e}_{\e s}$, 
  $\hat{Y}^{\e}_s=Y^{\e}_{\e s}$,
   %$\hat{Z}^{\e}_s=Z^{\e}_{\e s}$,
    $\hat{\Xi}^{\e}_s=\e^{-1/2}\Xi^{\e}_{\e s}$
    then  the last two equations in \eqref{SystemEpsilon} becomes:
 \begin{equation}\label{SystemEpsilonSlow}
\begin{cases}
%dX_{\e s}=  \e AX_{\e s} \, ds+ \sqrt{\e} R  \, d\hat{W}^1_s, & \\
  d\hat{Q}^\e_s= (B\hat{Q}^\e_s+ F(X_{\e s},\hat{Q}^\e_s) \,dt +   \,G \, d\hat{W}^2_s, &  \\ 
- d\hat{Y}^{\e}_s  =  \psi(X_{\e s},\hat{Q}^\e_s,{Z}^\e_{\e s}, \hat{\Xi}^\e_s) \,ds - \sqrt{\e} {Z}^\e_{\e s} d\, \hat{W}^1_s  - \hat{\Xi}^\e_s d\hat{W}_s^2, \\ 
X_0= x_0 \quad \hat{Q}^\e_0=q_0, \quad \hat{Y}^\e_{1/\e}=h(X_1).
\end{cases} 
\end{equation}
where  $\hat{W}_s^{\ell}={\e}^{-1/2}{W}_{\e s }^{\ell}$, $\ell=1,2$. We will often make use of this change of time in the proof.

We must compare:
$$
Y^\e _0 -\bar Y_0 =+ \int_0^1(\psi(X_t,Q^\e_t,Z^\e_t, \Xi^\e_t/\sqrt{\e} ) - \lambda({X}_t,\bar{Z}_t))\, dt -\int_0^1 (Z^\e_t-\bar{Z}_t)\,d W^1_t  -  \int_0^1\Xi^\e_t\,d W^2_t .$$

By adding and subtracting we split the first integral on the right hand side as:
\begin{align}\nonumber
 \int_0^1  (\psi(X_t,Q^\e_t,Z^\e_t,\Xi^\e_t/\sqrt{\e})-\lambda({X}_t,\bar{Z}_t))\,dt & =  
 \int_0^1  [(\psi({X}_t,Q^\e_t,{Z}^\e_t,\Xi^\e_t/\sqrt{\e})-\psi({X}_t,Q^\e_t,\bar{Z}_t,\Xi^\e_t/\sqrt{\e})] \,dt
\\ \label{primo_split}&+\int_0^1  (\psi({X}_t,Q^\e_t,\bar{Z}_t,\Xi^\e_t/\sqrt{\e}) -\lambda({X}_t,\bar{Z}_t)\,dt. \end{align}

We have to use a discretization argument to cope with the second member of the sum.

%{ \bf   $I_2$ term.}

Let us now introduce for every $ N$ positive integer, a partition of the interval $[0,1]$ of the form $t_k= {k}{2^{- N}}, \, k=0,1, \dots, 2^N$ and define a couple of  step processes ${X}^N$ and $\widetilde{Z}^N$ defined as follows:
\begin{equation}\label{defXN}
{X}^N(t)= {X}( t_k),\quad {t} \in[t_k, t_{k+1}[, \quad k =0,\dots 2^N-1, \end{equation}
\begin{equation}\label{defZtilde}
\widetilde{Z}^N(t)= 2^{N} \int_{t_{k-1}}^{t_k}\bar{Z}_{\ell} \, d\ell,  \text{ for } t  \in[t_k, t_{k+1}[, \ k =1,\dots 2^N-1, \ \widetilde{Z}(t)= 0 \ \text{ for } t \in[0, t_{1}[,
\end{equation}
where $X, \bar{Z}$ are part of the solution of \eqref{LimitEquation}. By construction one has that:
\begin{equation}\label{RegBarZ}
\lim_{N \to \infty}\E \int_{0}^{1} |\widetilde{Z}^N_t-\bar{Z}_t|^2 \, dt =0.
\end{equation}
We fix $N$, then for  $ k=0,1, \dots ,2^{N}-1$ we consider the following, iteratively defined, class  of forward SDE:
\begin{equation} \label{quickN}
d{\hat{\mathcal{Q}}}^{N,k}_s=(B{\hat{\mathcal{Q}}}^{N,k}_s+ F({X}_{  t_k} ,{\hat{\mathcal{Q}}}^{N,k}_s))\,dt + G d\hat{W}^2_s, \quad s \geq t_k/\e,
\hspace{20pt} \hat{\mathcal{Q}}^{N,k}_{t_k/\e}={\hat{\mathcal{Q}}^{N,k-1}_{t_k/\e}},
\end{equation}
Moreover we define (see Theorem \ref{Teolambda}):
 $$\check{Y}^{N,k}_s=\check{v}(X_{t_k},\hat{\mathcal{Q}}^{N,k}_s, \widetilde{Z}^N_{t_k}), \quad 
 \check{\Xi}^{N,k}_s=\check{\zeta}(X_{t_k},\hat{\mathcal{Q}}^{N,k}_s, \widetilde{Z}^N_{t_k}), \quad \hbox{ for } s\geq t_k/\e ,$$
so that
% by Markovianity of the process $\hat{\mathcal{Q}}^{N,k}$ 
the triplet $((\check{Y}^{N,k}_s)_{s\geq t_{k/\e}}, \lambda(X_{t_k},\widetilde{Z}^N_{t_k}), (\check{\Xi}^{N,k}_s)_{s\geq t_{k/\e}})$ verifies:
\begin{equation}
\label{BackInfHorStep}
- d \check{Y}^{N,k}_s=[\psi({X}_{t_k},\hat{\mathcal{Q}}^{N,k}_s,\widetilde{Z}^{N}_{t_k},\check{\Xi}^{N,k}_s)-\lambda({X}_{t_k},\widetilde{Z}^{N}_{ t_k})]\,\rmd s - \check{\Xi}^{N,k}_s \rmd \hat{W}_s^{2}, \quad  \text{ for all } \, s\geq t_{k}/\e,
\end{equation}
and
\begin{equation}\label{stimapathwise}
|\check{Y}^{N,k}_s|\leq c( 1 + |\widetilde{Z}^N_{t_k}|)
|\hat{\mathcal{Q}}_s^{N,k}|, \qquad \text{ for all } s  \geq t_k/\e, \end{equation}
for some positive constant $c>0$ independent of $k$ and $N$.

We also set for $  s\in [0,1/\e[$:
\begin{equation}\label{def-QNXIN}
\hat{\mathcal{Q}}^N_s=\sum_{k=0}^{2^N-1} \hat{\mathcal{Q}}^{N,k}_s I_{[t_k/\e, t_{k+1}/\e[}(s), \; \qquad\check{\Xi}^N_s=\sum_{k=0}^{2^N-1} \check{\Xi}^{N,k}_s I_{[t_k/\e, t_{k+1}/\e[}(s) ,\end{equation}
so that, for all $N\in \mathbb{N}$ and $k=0,..., 2^N-1$ have:
\begin{align}\label{diffYuguale0}
\check{Y}^{N,k}_{t_{k}/\e}-\check{Y}^{N,k}_{t_{k+1}/\e }-\int_{t_{k}/\e }^{t_{k+1}/\e}[\psi({X}^N_{\e s},\hat{\mathcal{Q}}^N_s,\widetilde{Z}^N_{\e s},\check{\Xi}^{N}_s)-\lambda({X}^N_{\e s},\widetilde{Z}^N_{\e s})]\, ds + 
\int_{t_{k}/\e }^{t_{k+1}/\e }\check{\Xi}^{N}_s \, d\hat{W}^2_s=0.
\end{align}
%therefore:
%\begin{align}
%\int_{0}^{1/\e}[\psi({X}^N_{\e t},\hat{Q}^N_t,\widetilde{Z}^N_{\e t},\hat{\Xi}^{N}_t)-\lambda({X}^N_{\e t},\widetilde{Z}^N_{ \e t})]\, dt = \sum_{k=1}^N \left[( \hat{Y}^{N,k}_{t_{k}/\e}-\hat{Y}^{N,k}_{t_{k+1}/\e })+
%\int_{t_{k} /\e}^{t_{k+1} /\e}\hat{\Xi}^{N}_t \, d\hat{W}^2_t\right].
%\end{align}
The second integral in the right hand side of \eqref{primo_split} can be written as:
$$
\int_0^1  (\psi({X}_t,Q^\e_t,\bar{Z}_t,\Xi^\e_t/\sqrt{\e}) -\lambda({X}_t,\bar{Z}_t)\,dt=\e \sum _{k=0}^{2^N-1}\int_{t_k/\e}^{t_{k+1}/\e}  [\psi(X_{\e s},\hat{Q}^\e_{ s},\bar{Z}_{\e s},\hat{ \Xi}^\e_{s}) -\lambda({X}_{\e s},\bar{Z}_{\e s})]\,ds,$$
and, adding the null terms in (\ref{diffYuguale0}) for $k=1,..., 2^N$, as:

\begin{align}\label{decI2}
&\int_0^1  (\psi({X}_t,Q^\e_t,\bar{Z}_t,\Xi^\e_t/\sqrt{\e}) -\lambda({X}_t,\bar{Z}_t)\,dt=\e \sum _{k=0}^{2^N-1}\int_{t_k/\e}^{t_{k+1}/\e}  [\psi(X_{\e s},\hat{Q}^\e_{ s},\bar{Z}_{\e s},\hat{ \Xi}^\e_{s}) - \psi({X}^N_{\e s},\hat{\mathcal{Q}}^{N}_{ s},\widetilde{Z}^N_{\e s},\check{\Xi}^N_s)]\,ds  \\&+\e \sum _{k=0}^{2^N-1}\int_{t_{k}/\e }^{t_{k+1}/\e }\check{\Xi}^{N}_s \, d\hat{W}^2_s 
-\e\sum _{k=0}^{2^N-1}\int_{t_k/\e}^{t_{k+1}/\e} [\lambda({X}_{\e s},\bar{Z}_{\e s})- \lambda(X^N_{\e s},\widetilde{Z}^N_{\e s})]\,ds \nonumber +\e  \sum_{k=1}^{2^N-1} (\check{Y}^{N,k}_{t_{k}/\e}-\check{Y}^{N,k}_{t_{k+1}/\e })
.\end{align}
Therefore coming back to our original term $Y^{\e}_0-\bar{Y}_0$ we have, taking into account \eqref{primo_split}:
\begin{align*}
Y^{\e}_0-\bar{Y}_0= & \, \e  \sum_{k=1}^{2^N-1} ( \check{Y}^{N,k}_{t_k/\e}-\check{Y}^{N,k}_{t_{k+1}/\e} )+\e\int_0^{1/\e} \left[\psi(X_{\e s }, \hat{Q}^{\e}_s, Z^{\e}_{\e s }, \hat{\Xi}^{\e}_s)-\psi(X_{\e s }, \hat{Q}^{\e}_s, \bar{Z}_{\e s }, \hat{\Xi}^{\e}_s)\right] \,ds\\&  
+\e \int_{0}^{1/\e}  [\psi(X_{\e s},\hat{Q}^\e_{ s},\bar{Z}_{\e s},\hat{ \Xi}^\e_{s}) - \psi({X}^N_{\e s},\hat{\mathcal{Q}}^{N}_{ s},\widetilde{Z}^N_{\e s},\check{\Xi}^N_s)]\,ds 
-\e\int_{0}^{1/\e} [\lambda({X}_{\e s},\bar{Z}_{\e s})- \lambda(X^N_{\e s},\widetilde{Z}^N_{\e s})]\,ds \\ &
-\sqrt{\e}\int_0^{1/\e } (Z^{\e}_{\e s}-\bar{Z}_{\e s}) d \hat{W}^1_s - {\e}\int_0^{1/\e } (\hat{\Xi}^{\e}_{ s}-\check{\Xi}^N_{ s}) d \hat{W}^2_s .
\end{align*}
Notice that we can rewrite this difference as follows:
\begin{equation}\label{secondo_split}
\begin{split}
 Y^\e _0 -\bar Y_0 = \, & 
 \e \int_0^{1/\e} {\mathcal{R}}^{\e,N}_s \, ds
 +\e  \sum_{k=1}^{2^N-1} (\check{Y}^{N,k}_{t_{k}/\e}-\check{Y}^{N,k}_{t_{k+1}/\e })\\ &     +\e \int_0^{1/\e}[\psi({X}_{\e t},\hat{{Q}}^{\e}_{s},{Z}^{\e}_{\e s},\hat{\Xi}^\e_s)-\psi({X}_{\e s},\hat{Q}^{\e}_{s},\bar{Z}_{\e s},\hat{\Xi}^\e_s)]\, ds \\ &  
  +\e \int_0^{1/\e}[\psi({X}^N_{\e s},\hat{\mathcal{Q}}^{N}_{s},\widetilde{Z}^N_{\e s},\hat{\Xi}^\e_s)-\psi({X}^N_{\e s},\hat{\mathcal{Q}}^N_s,\widetilde{Z}^N_{\e s},\check{\Xi}^N_t)]\, ds
 \\ & 
  +  \e \int_0^{1/\e}(\check{\Xi}^{N}_s -\hat{\Xi}^\e_s)  \, d\hat{W}^2_s +\sqrt{\e}\int_0^{1/\e} (Z^\e_{\e s}-\bar{Z}_{ \e s})\,d \hat{W}^1_s,
\end{split}
\end{equation}
where  $ {\mathcal{R}}^{\e,N}_s := \psi({X}_{\e s},\hat{Q}^{\e}_{s},\bar{Z}_{\e s},\hat{\Xi}^\e_s) - \psi({X}^N_{\e s},\hat{\mathcal{Q}}^{N}_{s},\widetilde{Z}^N_{\e s},\hat{\Xi}^\e_s)$.
Then by  Hypothesis \ref{B.3}   we deduce that for a suitable constant $c$, independent from $\e$ and $N$, the following holds:
\begin{equation}\label{primoterm}
|\mathcal{R}^{\e,N}_s|\leq c (1+|\bar{Z}_{\e s}|)|X_{\e s }-X^N_{\e s}|+c(1+|\bar{Z}_{\e s}|)|\hat{Q}^{\e}_s-\hat{\mathcal{Q}}^N_s|+ c|\bar{Z}_{\e s }-\widetilde{Z}^N_{\e s}|.
\end{equation}

\medskip

\noindent The presence of the two stochastic in \eqref{secondo_split} allows us to get rid of the third and fourth term on the right hand side by  a Girsanov argument, namely we introduce:
\begin{equation}\label{Girsa1}
\delta^{1,\e} (s) = \begin{cases}\displaystyle\frac{[\psi({X}_{\e t},\hat{Q}^{\e}_{s},{Z}^{\e}_{\e s},\hat{\Xi}^\e_s)-\psi({X}_{\e s},\hat{Q}^{\e}_{s},\bar{Z}_{\e s},\hat{\Xi}^\e_s)]}{|{Z}^\e_{\e s} -\bar{Z}_{ \e s}|^2} ({Z}^\e_{\e s} -\bar{Z}_{ \e s})^*& \text { if } |{Z}^\e_{\e s} -\bar{Z}_{\e s}| \not= 0, \\
0 & \text { if } |{Z}^\e_{\e s} -\bar{Z}_{ \e s}| = 0,
 \end{cases}
\end{equation}
and 
\begin{equation}\label{Girsa2}
\delta^{2,\e,N} (s) = \begin{cases}\displaystyle\frac{\psi({X}^N_{\e s},\hat{\mathcal{Q}}^{N}_{s},\widetilde{Z}^N_{\e s},\hat{\Xi}^\e_s)-\psi({X}^N_{\e s},\hat{\mathcal{Q}}^N_s,\widetilde{Z}^N_{\e s},\check{\Xi}^N_s)}{|\hat{\Xi}^\e_s -\check{\Xi}^N_s|^2}(\hat{\Xi}^\e_s -\check{\Xi}^N_s)^* & \text { if } |\hat{\Xi}^\e_s -\check{\Xi}^N_s| \not= 0, \\
0 & \text { if } |\hat{\Xi}^\e_s -\check{\Xi}^N_s|= 0.
 \end{cases}
\end{equation}
We notice that processes $(\delta^{1,\e}(s))_{s\in [0,1/\e]}$ and 
$(\delta^{2,\e, N}(s))_{s\in [0,1/\e]}$ are  bounded uniformly by $L_\xi$ and $L_z$ respectively, see Hypothesis \ref{B.3}. 
We have:
\begin{align*}
Y^\e _0 -\bar Y_0 = & \, \e \int_0^{1/\e} \delta^{1,\e}( s) [Z^\e_{\e s} -\bar{Z}_{\e s}] \, ds +\e \int_0^{1/\e}\delta^{2,\e,N}(s)[\check{\Xi}^{N}_s -\hat{\Xi}^\e_s]\, ds  \\ &
 +  \e \int_0^{1/\e}(\check{\Xi}^{N}_s -\check{\Xi}^\e_s)  \, d{\hat{W}}^2_s  +\sqrt{\e}\int_0^{1/\e} (Z^\e_{\e s}-\bar{Z}_{ \e s})\,d \hat{W}^1_s
\\ & 
 +\e \int_0^{1/\e} {\mathcal{R}}^{\e,N}_s \, ds +\e  \sum_{k=1}^{2^N-1} (\check{Y}^{N,k}_{t_{k}/\e}-\check{Y}^{N,k}_{t_{k+1}/\e }).
\end{align*}
and rescaling time (speeding it up this time)
\begin{align*}
Y^\e _0 -\bar Y_0 = & \, \int_0^{1} \delta^{1,\e}( t/\e) [Z^\e_{ t} -\bar{Z}_{ t}] \, dt + \int_0^{1}\delta^{2,\e,N}(t/\e)[\check{\Xi}^{N}_{\e^{-1}t} -\hat{\Xi}^\e_{\e^{-1}t}]\, dt 
\\ &
 +  \sqrt{\e} \int_0^{1}(\check{\Xi}^{N}_{\e^{-1}t} -\check{\Xi}^\e_{\e^{-1}t})  \, d{{W}}^2_t  +\int_0^{1} (Z^\e_{ t}-\bar{Z}_{ t})\,d W^1_t
\\ & 
 + \int_0^{1} {\mathcal{R}}^{\e,N}_{\e^{-1}t} \, dt +\e  \sum_{k=1}^{2^N-1} (\check{Y}^{N,k}_{t_{k}/\e}-\check{Y}^{N,k}_{t_{k+1}/\e }) .
\end{align*}
We set, for $ t\in [0,1]$:
\begin{align}
&\widetilde{W}^{1}_t =:  \int_0^{t} \delta^{1,\e}(r/\e)\,dr +  W^1_t, \\
&\widetilde{W}^2_{t}=:  \e^{-1/2}\int_0^{t}\delta^{2,\e,N}(r/\e)\,dr +  W^2_{t}.
\end{align}
We denote by $\widetilde{\E}^\e$ the expectation under the new probability $\tilde{\mathbb{P}}^\e$  with respect to which $(\widetilde{W}^{1}_t, \widetilde{W}^2_{t})_{t\in [0,1]}$ is a $H\times K$ valued cylindrical Wiener process (recall that $({W}^{1}_t, {W}^2_{t})_{t\in [0,1]}$ is a $H\times K$ valued cylindrical Wiener process). 
Since the left hand side is deterministic, we have:
\begin{equation}\label{diffGirs}
Y^\e _0 -\bar Y_0 = \widetilde{\E}^\e \int_0^{1} {\mathcal{R}}^{\e,N}_{t/\e} \, dt+ \e\widetilde{\E}^\e  \sum_{k=1}^{2^N-1} [\check{Y}^{N,k}_{t_{k}/\e}-\check{Y}^{N,k}_{t_{k+1}/\e }].\end{equation}
Moreover, taking into account \eqref{primoterm}, it holds: $$ \widetilde{\E}^\e \int_0^{1} |{\mathcal{R}}^{\e,N}_{t/\e}| \, dt \leq c \widetilde{\E}^\e \int_0^{1}  \left((1+|\bar{Z}_t|)|{X}_{ t}-{X}^N_t| + 
(1+|\bar{Z}_t|)|\hat{Q}^\e_{ t/\e}-\hat{\mathcal{Q}}^N_{t /\e}|+ |\bar{Z}_{t}-\widetilde{Z}^N_t|\right)\, dt.$$

Let us start from 
 $$\widetilde{\E}^\e \int_0^{1}  (1+|\bar{Z}_t|)|{X}_{ t}-{X}^N_t|dt.$$
 We notice that, with respect to $\widetilde{W}^1$ we have:
 \begin{equation*}
\begin{cases}
d X_t =   AX_t \, dt - R \delta^{1,\e}({t}/{\e})dt + R\, d \widetilde{W}^1 _t, \\
-d\bar{Y}_t  = \lambda({X}_t,\bar{Z}_t)\,dt - \, \bar{Z}_t[- \delta^{1,\e}({t}/{\e})dt+\, d\widetilde{W}^1_t],   \\
\bar{Y}_1= h(X_1), \quad X_0= x_0.
\end{cases} 
\end{equation*}
% $(X_t)_{t\geq 0}$ satisfies 
%$$dX_t= AX_tdt -\delta^{1,\e}({t}/{\e})dt+ \, d%\widetilde{W}^1_t, \quad X_0=x^0 $$
Define: $$\rho:=\exp\left( \int_0^1 \delta^{1,\e}(s/\e) \, d \widetilde{W}_s^{1}  - \frac{1}{2} \int_0^1 | \delta^{1,\e}(s/\e) |^2 \, ds \right),$$
then, by Holder inequality, setting $\Delta_{X,N}:=
\sup_{t\in [0,1]}|{X}_{ t}-{X}^N_t|$ it holds:
\begin{align*}
\widetilde{\E}^\e \int_0^{1} (1+|\bar{Z}_t|) |{X}_{ t}-{X}^N_t|dt\leq \widetilde{\E}^\e\left[ \Delta_{X,N} \int_0^{1} (1+|\bar{Z}_t|) dt \right] \qquad\qquad\qquad\qquad \qquad\qquad \\\leq
\nonumber\widetilde{\E}^\e\left[ \rho^{-3/4}( \rho^{1/4}\Delta_{X,N}) \rho^{1/2}\!\!\!\int_0^{1}(1+|\bar{Z}_t|)dt\right]
\leq
\left[ \widetilde{\E}^\e \rho^{-3}\right]^{1/4} \left[ \widetilde{\E}^\e (\rho  \Delta_{X,N}^4)\right]^{1/4}\left[ \widetilde{\E}^\e \left(\!\rho\! \int_0^1\! (1+|\bar{Z}_t|^2 )dt\right) \right]^{1/2}.
\end{align*}

Again by Girsanov the process
$\left(-\int_0^t \delta^1(t/\e) dt+ \widetilde{ W}^1_t\right)_{t\in [0,1]}$ is a cylindrical Wiener process with respect to $\rho\, d\mathbb{P}^{\e}$. By uniqueness of the solution of the forward backward system (\ref{LimitEquation}) the law of the process $(X_t)_{t\geq 0}$  under $\rho d\mathbb{P}^{\e}$  coincides with its law with respect to $\mathbb{P}$.  Moreover  we notice that being $\bar{Z}_t= \zeta (X_t)$ where $\zeta$ is a deterministic Borel function $H\rightarrow \Xi^*$ then the law of $\bar{Z}$ and $\tilde{Z}^N$ depend only on the law of $(X)$ in a non anticipating way.  So even the law of $(\bar{Z}_t)_{t\geq 0}$ and  $(\bar{Z}^N_t)_{t\geq 0}$ under $\rho d\mathbb{P}^{\e}$  coincides with its law with respect to $\mathbb{P}$. 

Recalling that $\delta^{1,\e}$ is uniformly bounded and consequently (with respect to $\e$ as well) we have $\widetilde{\E}^\e \rho^{-3}\leq c$ (where $c$ does not depend on $\e$), moreover   
$$\widetilde{\E}^\e\left(\rho\int_0^1 |\bar{Z}_t|^2 dt \right)={\E}\left(\int_0^1 |\bar{Z}_t|^2 dt \right)<+\infty.
$$
Thus we can conclude
\begin{equation} \label{stimaZNX}\widetilde{\E}^\e \int_0^{1} (1+|\bar{Z}_t|) |{X}_{ t}-{X}^N_t|\, dt\leq
C[\E \Delta_{X,N}^4]^{1/4},
\end{equation}
where $C$ is independent of $N$ and $\e$. 

\noindent By the continuity of trajectories of $(X_t)_{t\geq 0}$, having also $\E \sup_{t\in [0,1]} |X_t|^4<\infty$,  we get: 
\begin{equation}\label{Delta1}
\E \Delta_{X,N}^4\rightarrow 0, \qquad \text{ as }\qquad N\rightarrow \infty.\end{equation} 
We also have that:

%\noindent Concerning the term 
% $$\widetilde{\E}^\e \int_0^{1}  |\bar{Z}_{ t}-\widetilde{Z}^N_t|dt$$
% we notice that being $\bar{Z}_t= \zeta (X_t)$ where $\zeta$ is a deterministic Borel function $H\rightarrow \Xi^*$ then the law of $(\bar{Z})$ and $(\tilde{Z}^N)$ depend only on the law of $(X)$ in a non anticipating way. So, the above  argument still works and yields:
 \begin{equation}\label{stimabarZ}
\widetilde{\E}^\e\int_0^{1}  |\bar{Z}_{ t}-\widetilde{Z}^N_t|dt
%={\E} \, \left[\rho^{-1}\!\!\!\int_0^{1}  |\bar{Z}_{ %t}-\widetilde{Z}^N_t|dt \right]
\leq
C  \left[\E \int_0^{1} |\bar{Z}_{ t}-\widetilde{Z}^N_t|^2dt\right]^{1/2}= C (\E\Delta_{Z,N})^{1/2},
\end{equation}
where $\Delta_{Z,N}= \displaystyle \int_0^{1} |\bar{Z}_{ t}-\widetilde{Z}^N_t|^2dt$ and by  \eqref{RegBarZ}:
\begin{equation}\label{Delta2}
\E \Delta_{Z,N}\to 0, \qquad \text{ as }\qquad N\rightarrow \infty.\end{equation} 
Now we deal with  the term:
\begin{align*}
 \widetilde{\E}^{\e}  \, \int_0^{1} (1+|\bar{Z}_t|) |\hat{Q}^\e_{\e^{-1} t}-\hat{\mathcal{Q}}^N_{ \e^{-1} t}|\, dt .
\end{align*}
%With respect to $\widetilde W^2 $ the process $({Q}^{\e}_t)_{t\in [0,1]}$ solves
%$$
%\e\, d{{Q}^{\e}}_t=(B{{Q}^{\e}}_t+ F({X}_{ t},{{Q}^{\e}}_t))\,dt - \delta^{2,\e,N}(t/\e)\,dt+ \sqrt{\e}\,d\widetilde{W}^2_t , \quad t \geq 0,
%\hspace{20pt} {{Q}^{\e}}_0={q_0},
%$$

Introducing the $\mathbb{P}^{\e}$ Wiener process $\hat{\widetilde{W}}_s:= (\e)^{-1/2}\widetilde{W}_{\e s}$ we have that the process $(\hat{Q}^{\e}_s)_{s\in [0,1/\e]}$ solves:
\begin{equation} \label{quicktilde}
\, d\hat{Q}^{\e}_s=(B{\hat{Q}^{\e}}_s+ F({X}_{ \e s},{\hat{Q}^{\e}}_s))\,ds -\delta^{2,\e,N}(s)\,ds-  Gd\hat{\widetilde{W}}^2_s , \quad s \geq 0,
\hspace{20pt} {\hat{Q}^{\e}}_0={q_0},
\end{equation}
moreover $(\hat{\mathcal{Q}}^{N}_{ s})_{s\in [0,1/\e]}$ solves:
\begin{equation} \label{quicktildeN}
d{\hat{\mathcal{Q}}^{N}}_s=(B{\hat{\mathcal{Q}}^N}_s+ F({X}^N_{\e s },{\hat{\mathcal{Q}}^{N}}_s))\,dt - \delta^{2,\e,N}( s)\,ds+ G d\hat{\widetilde{W}}^2_s, \quad s \geq 0,
\hspace{20pt} \hat{\mathcal{Q}}^{N}_0={q_0}.
\end{equation}
    Therefore by Lemma \ref{dipendcont} and hypothesis  \ref{B.3} we have for all $p\geq 1$:
\begin{equation}
\sup_{s \in [0, 1/\e] } \widetilde{\E}^\e[|\hat{\mathcal{Q}}^{N}_t|^p]\leq c_p\left(1+|q_0|^p+\sup_{s \in [0, 1/\e]} \widetilde{\E}^\e |X_s|^p + \sup_{s \in [0, 1/\e]} \widetilde{\E}^\e \left | \int_0^{s} e^{(s-r) B}   G d\hat{\widetilde{W}}^2_r \right|^p+ L_{\xi}\right),
\end{equation}
for a constant $c_p$ independent of  $\e$ and $N$.
Arguing as before, we have that
\[  \widetilde{\E}^\e |X_s|^p= \widetilde{\E}^\e ( \rho^{-1/2}  \rho^{1/2}|X_s|^p) \leq   (\widetilde{\E}^\e  \rho ^{-1})^{1/2}    (\widetilde{\E}^\e  ( \rho|X_s|^{2p}))^{1/2} \leq C 
  (\E  |X_s|^{2p})^{1/2}, \]
 and 
 \begin{align*} 
 & \displaystyle\widetilde{\E}^\e\left|\int_0^{s} e^{(s-r) B}   G d\hat{\widetilde{W}}^2_r\right|^p= \widetilde{\E}^\e \left( \rho^{-1/2}  \rho^{1/2}\left|\int_0^{s} e^{(s-r) B}   G d\hat{\widetilde{W}}^2_r\right|^p \right)  \\
  & \leq  (\widetilde{\E}^\e  \rho ^{-1})^{1/2}    \left(\widetilde{\E}^\e  \left( \rho\left|\int_0^{s} e^{(s-r) B}   G d\hat{\widetilde{W}}^2_r\right|^{2p}\right)\right)^{1/2} \leq C 
  \left(\E  \left|\int_0^{s} e^{(s-r) B}   G d\hat{\widetilde{W}}^2_r\right|^{2p}\right)^{1/2}, 
  \end{align*}
  for some constant $C>0$ independent of $\e$.
Therefore,   bearing in mind the estimate \eqref{stimaslow} for the slow component $X$ and hypothesis  \ref{A.5},  we conclude that there exists a constant $c >0$, independent of $\e $ and $N$, such that
\begin{equation}\label{stimaunifQ}
\sup_{s \in [0, 1/\e] } \widetilde{\E}^\e[|\hat{\mathcal{Q}}^{N}_t|^p]\leq c.
\end{equation}
Again by  Lemma \ref{dipendcont} one has that for all $s>0$, $$|\hat{{Q}}^\e_{ s}-\hat{\mathcal{Q}^{N}_{ s}}|\leq c\int_0^s e^{-\eta(s-\ell)}|X_{\e \ell}-X^N_{\e \ell}| d \ell\leq c\Delta_{X,N}, $$
thus, arguing as in \eqref{stimaZNX}, 
\begin{equation}\label{stimaZQ}
\displaystyle\widetilde{\E}^\e \int_0^{1}(1+|\bar{Z}_t|) |\hat{{Q}}^\e_{ \e^{-1}t}-\hat{\mathcal{Q}}^{N}_{ \e^{-1}t}|dt\leq c \widetilde{\E}^\e\left[ \Delta_{X,N}\int_0^1(1+|\bar{Z}_t|) d t\right]\leq C[\E \Delta_{X,N}^4]^{1/4},
\end{equation}
as above.

Now we come to the last term.
%\begin{align*}
%\e \widetilde{\E}^\e  \sum_{k=1}^N ( \hat{Y}^{N,k}%_{t_{k}/\e}-\hat{Y}^{N,k}_{t_{k+1} /\e})= \e %\widetilde{\E} ^\e \sum_{k=1}^N \hat{Y}^{N,k}_{t_{k}/%\e}-\hat{Y}^{N,k}_{t_{k+1} /\e}| \F_{t_{k}})
%\end{align*}
%Being $\mathbb{Q}^\e$ absolutely continuous with %respect to $P$, \eqref{stimapathwise} guarantees:
%\begin{equation}
%\widetilde{\E}^{\F_{t_k}}(|\hat{Y}^{N,k}_t|)\leq c( 1 %+\widetilde{\E}^{\F_{t_k}} |\hat{Q}_t^{N}|) \qquad %%\text{ for all } t  \geq t_k \quad \mathbb{Q}^\e-%a.s.
%\end{equation}
%Notice that, under $\mathbb{Q}^\e$ the process  $\hat{Y}^{N,k}$ has the following differential
%\begin{align*}
%\hat{Y}^{N,k}_{t_{k+1}/\e}-\hat{Y}^{N,k}_{t_{k}/\e }=\int_{t_{k}/\e }^{t_{k+1}/\e}[\psi(X^N_{t_k},Q^N_t, Z^{N}_{t_k},\hat{\Xi}^{N}_t)-\lambda(X^N_{t_k}, Z^N_{t_k})+ \delta^2\psi^{N,\e}(t)\hat{\Xi}^{N}_t]\, dt + 
%\int_{t_{k}/\e }^{t_{k+1}/\e }\hat{\Xi}^{N}_t \, d\widetilde{W}^2_t.
%\end{align*}
%We have that:
%\begin{equation*}
%|\widetilde{\E}^{\F_{t_k}} ( \hat{Y}^{N,k}_{t_{k}/\e}-%\hat{Y}^{N,k}_{t_{k+1} /\e})| \leq \widetilde{\E}%^{\F_{t_k}}  ( |\hat{Y}^{N,k}_{t_{k}/\e}|+|\hat{Y}%^{N,k}_{t_{k+1} /\e}|)\leq C (1+\widetilde{\E}%^{\F_{t_k}} |\hat{Q}^N_{t_k/\e}|+\widetilde{\E}%^{\F_{t_k}} \hat{Q}^N_{t_{k+1}/\e}|)
%\end{equation*}
    We apply successively (\ref{stimapathwise})  and \eqref{stimaunifQ} to get the following estimates (the value of the constant $c$ below can change from line to line but never depends neither on $k$ nor on $N$ or on $\e$):
\begin{align*}
&\left|\e \widetilde{\E}^\e  \sum_{k=1}^{2^N-1} ( \hat{Y}^{N,k}_{t_{k}/\e}-\hat{Y}^{N,k}_{t_{k+1} /\e})\right|\leq c \e  \sum_{k=1}^{2^N-1} \widetilde{\E}^{\e}\left[ (1+|\widetilde{Z}^N_{t_k}|)(1+ |\hat{\mathcal{Q}}^N_{t_k/\e}|+|\hat{\mathcal{Q}}^N_{t_{k+1}/\e}|)\right]\\
& \leq  c\e \sum_{k=1}^{2^N-1} \left[\widetilde{\E}^{\e}  (1+|\widetilde{Z}^N_{t_k}|)^{4/3}\right]^{3/4}\left[\widetilde{\E}^{\e}(1+ |\hat{\mathcal{Q}}^N_{t_k/\e}|+|\hat{\mathcal{Q}}^N_{t_{k+1}/\e}|)^4\right]^{1/4}\leq c \e \sum_{k=1}^{2^N-1} \left[1+\left(\widetilde{\E}^{\e}  |\widetilde{Z}^N_{t_k})^{4/3}\right)^{3/4}\right].
\end{align*}
Proceeding as above, recalling that the law of  $\widetilde{Z}^N_{t_k}$ depends only on the law of the process $(X_t)$ we have:
$$\widetilde{\E}^{\e} [(|\widetilde{Z}^N_{t_k}| )^{4/3} ]\leq \left[{\E} \rho^{-2}\right]^{1/3} \left[\E |\widetilde{Z}^N_{t_k}|^{2}\right]^{2/3} \leq   c 2^{\frac{ 2}{3}N}\left [ \E\int_{0}^{t} |\bar{Z}_{t}|^2 dt  \right]^{2/3}.$$
At last we sum up the latter result,  \eqref{stimaZNX}, \eqref{stimabarZ} and \eqref{stimaZQ} to get:
\begin{align*}
&|Y^\e _0 -\bar Y_0 |\leq \widetilde{\E}^\e \int_0^{1} |\mathcal{R}^{\e,N}_{t/\e}| \, dt+ \e\widetilde{\E}^\e  \sum_{k=1}^{2^N}  |\hat{Y}^{N,k}_{t_{k}/\e}-\hat{Y}^{N,k}_{t_{k+1}/\e }|\\& \leq C[\E \Delta_{X,N}^4]^{1/4}+ C (\E\Delta_{Z,N})^{1/2}+ \e c  2^{\frac{3}{2} N}\left(\E\int_0^1 |\bar{Z}_t|^2 \, dt \right)^{1/2} +\e c 2^N.
\end{align*}
So letting first $\e$ tend to $0$ and then $N$ to $\infty$ the claim follows, by \eqref{Delta1} and \eqref{Delta2}.
\finedim
%%%%%%%%
\begin{remark} {\em Consider the following class of forward backward systems with initial time $\tau\in [0,1]$
 \begin{equation} \label{LimitEquationinitialdata}
\begin{cases}
d X^{\tau,x}_t =   AX^{\tau,x}_t \, dt + \, R d W^1 _t,\quad t\leq 1,  \\
-d\bar{Y}^{\tau,x}_t  =\lambda({X}^{\tau,x}_t,\bar{Z}^{\tau,x}_t)\,dt - \, \bar{Z}^{\tau,x}\, dW^1_t,\quad t\leq 1,   \\
X^{\tau,x}_\tau= x, \quad \bar{Y}_1^{\tau,x}=h(X^{\tau,x}_1).
\end{cases} 
\end{equation}  
If we set $v(\tau,x) = \bar{Y}^{\tau,x}_\tau$ then it is shown in \cite{FuTes_BE} that $v$ is a deterministic continuous function $[0,1]\times H\rightarrow \mathbb{R}$ being Gateaux differentiable with respect to the second variable. Moreover it is the unique mild solution of the nonlinear Kolmogorov equation
$$ \left\{\begin{array}{l}\displaystyle
\frac{\partial v(t,x)}{\partial t}+\mathcal{L} v(t,x) = \lambda (
x,\nabla v(t,x)R),\qquad t\in [0,1],\,
x\in H,\\
\displaystyle v(1,x)=h(x),
\end{array}\right.
$$
where $\mathcal{L}$ is the second order operator $$\mathcal{L}g(x)=\frac{1}{2}\hbox{Tr}[R^*\nabla^2 g(x)R ],\quad g\in \mathcal{C}^2(H),$$
$\nabla^2 g(x) \in \mathcal{L}(H)$ being the second derivative of $g$ in $x$.

In particular the limit $\lim_{\e\rightarrow 0}Y^{\e}_0$ can also be represented by the solution of the above HJB  equation as:
$$\lim_{\e\rightarrow 0}Y^{\e}_0=\bar{Y}^{0,x_0}_0=v(0,x_0).$$
}\end{remark} 

\section{The two scale control problem}

In this section we are finally in a position to exploit the convergence of solutions of BSDEs proved in the previous sections to solve our original problem of characterizing the limit of value functions $
V^{\e}(x_0,q_0)$ as $\e\rightarrow 0$ (see \eqref{intro-vf} in the Introduction).

\medskip

It turns out to be convenient to formulate the control problems in a weak form. Remark \ref{rem-equiv-weak} below reminds the reader about about the relation with the original formulation.

\medskip

Given the solution $(X,Q^{\e})$ of system \eqref{eqX} and a progressive measurable process $(\alpha_t)_{t\in [0,1 ]}$ taking its values in a complete metric space $U$ we denote by $\Theta^{\e,\alpha}$ the density
$$\Theta^{\e,\alpha}=\exp\left(\int_0^1 \left[ R^{-1}b(X_t,Q^\e_t,\alpha_t) dW^1_t + \frac{1}{\sqrt{\e}}\rho(\alpha_t)dW^2_t\right]-\frac{1}{2}\int_0^1 \!\!
\left[|R^{-1}b(X_t,Q^\e_t,\alpha_t)|^2+\frac{1}{\e}|\rho(\alpha_t)|^2\right]dt\right),$$
where $b: H\times K\times U \rightarrow H$ and $\rho: U \rightarrow K$ are measurable functions satisfying suitable assumptions listed below.

We also consider the following cost functional:
\begin{equation} \label{cost_epsilon}
J^{\e}(x_0,q_0,\alpha)= \E\left[\Theta^{\e,\alpha}\left( \int_0^1 l(X_t, Q^\e_t, \alpha_t)dt + h (X_1)\right)\right],
\end{equation}
where $l: H\times K\times U\rightarrow \mathbb{R}$ and $h: H \rightarrow \mathbb{R}$ are measurable and  satisfy the assumptions below:
\begin{hypothesis}\label{C.1}
There are positive constants $L$ and $M$ such that :
$$ |b(x,q,u)- b(x',q',\alpha)| \leq L (|x-x'|+|q-q'|),\qquad \qquad \forall \, q,q' \in K,  x, x' \in H, \, \alpha \in U,$$
$$ |l(x,q,\alpha)- l(x',q',\alpha)| \leq L (|x-x'|+|q-q'|),\qquad \qquad \forall \, q,q' \in K,  x, x' \in H, \, \alpha \in U,$$
$$ |h(x)- h(x')| \leq L |x-x'|,\qquad \qquad \forall x, x' \in H, $$
$$|b(x,q,\alpha)|, |l(x,q,\alpha)|, |\rho(\alpha)|, |h(x)| \leq M,\qquad \qquad \qquad  \forall q \in K, x \in H, \, \alpha \in U. $$
\end{hypothesis}
\begin{remark} \label{rem-equiv-weak}
We recall that if $d\mathbb{P}^{\e,\alpha}:=\Theta^{\e,\alpha} d\mathbb{P}$ then under probability $\mathbb{P}^{\e,\alpha}$ the process:
$$(\mathcal{W}^1_t,\mathcal{W}^2_t)=
(-\int_0^t R^{-1}b(X_r,Q^{\e}_r,\alpha_r)dr+{W}^1_t,-  \frac{1}{\sqrt{\e}}\int_0^t \rho(\alpha_r)dr+W^2_t),$$
is a cylindrical Wiener process in $\Xi\times \Xi$ and that with respect to $(\mathcal{W}^1_t,\mathcal{W}^2_t)$ the couple of processes $(X_t, Q^\e_t)$ satisfies the controlled system:
\begin{equation}
\label{eqcontr}
\begin{cases}
dX_t= AX_t \, dt  +  b(X_t,Q^{\e}_t,\alpha_t)dt \,+ R d\mathcal{W}^1_t, & X_0=x_0, \\ \\
 \e dQ^\e_t= (BQ^\e_t+ F(X^\e_t,Q^\e_t) \,dt +  G\rho(\alpha_t)dt+ \e^{1/2} G \, d \mathcal{W}^2_t, & Q^\e_0=q_0.
\end{cases}
\end{equation}
Moreover:
$$
J^{\e}(x_0,q_0,\alpha)= \E^{\mathbb{P}^{\e,\alpha}}\left( \int_0^1 l(X_t, Q^\e_t, \alpha_t)dt + h (X_1)\right),$$
thus the one introduced here is just a reformulation of our original control problem see
\eqref{controllointro} and \eqref{costointro}.
\end{remark}
We define, for $x\in H$, $q\in K$ and
$z, \xi \in \Xi^*$ :
\begin{equation}
\label{Hamiltonian}
\psi(x,q,z,\xi) = \inf_{\alpha \in U} \{l(x,q,\alpha)+ z[R^{-1}b(x,q,\alpha)] + \xi \rho(\alpha) \},  \end{equation}
and notice that, by straight forward cosiderations, under  Hypothesis \ref{C.1}, the Hamiltonian $\psi$ verifies hypothesis \ref{B.3}.
%%\begin{remark} For an eventual extension by approximation to the case of a degenerate $R$ we notice that if we define
% for $x\in H$, $q\in K$,
%$p\in H^*$, $\xi\in \Xi^*$ 
%\begin{equation}
%%\label{Hamiltonian-cancelletto}
%%\psi^{\sharp}(x,q,p,\xi) = \inf_{\alpha \in U} \{l(x,q,\alpha)+ pb(x,q,\alpha) + \xi \rho(v) \}  \end{equation}
%then $\psi(x,q,z,\xi)=\psi^{\sharp}(x,q,[R^{-1}]^*z,\xi) $
%where $[R^{-1}]^*: \Xi^*\rightarrow H^*$ is the adjoint of the right-inverse of $R$.
%\end{remark}

The main result of this paper is now just an immediate consequence of our general of Theorem \ref{main}.
\begin{theorem}\label{main_controllo} Denote by $V^{\e}$ the value function of our control problem, that is:
$$V^\e(x_0,q_0):= \inf_{\alpha}J^\e(x_0,q_0,\alpha),$$
where the infimum is taken over all progressive processes $\alpha$ with value in $U$. 

The sequence $V^\e(x_0,q_0)$ converges to the solution $\bar{Y}_0$ of equation \eqref{LimitEquation} evaluated at zero.
\end{theorem}

\Dim In \cite{FuHuGenGrad} it is shown that $V^\e(x_0,q_0)=Y^\e_0$ (see \eqref{SystemEpsilon}). The claim then follows by Theorem \ref{main}.

%\begin{hypothesis}\label{AssContr} 
%\item $|b(t,x,q)|\leq C_b,\qquad |b(t,x,q)-%b(t,x',q')|\leq L_b (|x-x'|+|q-q'|)$ for all $t\in %%[0,1]$, $x,x'\in H$, $q,q'\in K$.
%\begin{hypothesis}
%$$ $$

%\textbf{ Resta da scrivere la caratterizzazione di %lambda}
%$$ $$

%\begin{remark}\label{structure}{\em
%As we've already mentioned in the introduction we cannot cope  with the structure condition in the slow equation, that is we cannot deal with systems of the form:

%\begin{equation}
%\label{eqcontr_structure}
%\begin{cases}
%dX_t= AX_t \, dt  +  b(X_t,Q^{\e}_t)\,dt +R \mu(\alpha_t)\,dt \,+ R d\mathcal{W}^1_t, & X_0=x_0, \\ \\
% \e dQ^\e_t= (BQ^\e_t+ F(X^\e_t,Q^\e_t) \,dt +  G\rho(\alpha_t)dt+ \e^{1/2} G \, d \mathcal{W}^2_t, & Q^\e_0=q_0,
%\end{cases}
%\end{equation}
%where $R$ is allowed to be degenerate. This restriction seems related with the technique we have used. The presence of the fast component $Q$ in the slow equation even after the Girsanov transformation, would force the process $X$ to depend on $\e$ as well in  \eqref{SystemEpsilon}, and this would indeed affect the convergence  of the slow component.  
%}
%\end{remark}

\begin{remark} \label{rem-lambdacontrol}
{\em The nonlinearity $\lambda$ in the limit equation \eqref{LimitEquation} has itself a control theoretic interpretation. 
Namely, fixed $x\in H$ and $z\in \Xi^*$, let us consider the following ergodic control problem with {\em state equation}
\begin{equation}\label{ergfrozenstate}
d \hat{Q}^{\beta}_s=B\hat{Q}^{\beta}_sds+ F(x, \hat{Q}^{\beta}_s)\, d s + G\rho(\beta_s) d s +  G d\hat {W}^2_s,
\end{equation}
and {\em ergodic cost functional} 
\begin{equation}\label{ergfrozencost}
 \check{J}(x,z,\beta)=\liminf_{\delta \to 0} \E \, \delta \int_0^{\infty} e^{-\delta s} [z R^{-1}b(x, \hat{Q}^{\beta}_s,\beta_s)+ l(x,\hat{Q}^{\beta}_s,\beta_s)] d s.
\end{equation}
Then $\lambda(x,z)$ is the value function of the ergodic control problem the we have just described, that is:
$$\lambda(x,z)=\inf_{\beta}  \check{J}(x,z,\beta),$$
where the infimum is taken over all progressive processes ${\beta}: [0,\infty[\rightarrow U$.
\smallskip

\noindent
Notice that, in particular, being the infimum of linear functionals, the map $z\rightarrow \lambda(x,z)$ is concave.
\smallskip

\noindent Moreover notice that the result was proven in \cite{FuHuTes} with $\liminf$ replaced by $\limsup$ in the definition \eqref{ergfrozencost} of the ergodic cost nevertheless, as it can be easily verified, this substitution is inessential in the argument reported in  \cite{FuHuTes}
%\noindent Finally, again for further use, we prefer to introduce, for $p\in H^*$
%\begin{equation}\label{ergfrozencost_bis}
% \check{J}^{\sharp}(x,p,\beta)=\liminf_{\delta \to 0} \E \, \delta %\int_0^{\infty} e^{-\delta s} [p b(x, \hat{Q}^{\beta}_s,\beta_s)+ l(x,\hat{Q}%^{\beta}_s,\beta_s)] d s.
%\end{equation}
%and $\lambda^{\sharp}(x,p)=\inf_{\beta}  \check{J}_{}\sharp(x,p,\beta)$
%observing that $\lambda(x,z)=\lambda^{\sharp}(x,(R^{-1})^*z)$
}\end{remark}

\begin{example}\label{esempio}{\rm 
We provide a simple example
 to which our result apply.
Let us consider the following two scale system of classical controlled reaction diffusion  SPDEs in one space dimension driven by space time white noises see, for instance \cite{daza2} Section 11.2 or :
\begin{equation}\label{sistemaesempio}\left\{\begin{array}{l}
\displaystyle\frac{\partial}{\partial t} u^\e(t,x) = \frac{\partial^2}{\partial x^2} u^\e(t,x) +
 b(u^\e(t,x),v^\e(t,x), \alpha(t,x)) +
 \sigma(x)  \frac{\partial}{\partial t} {\mathcal{W}}^1(t,x),   \\ \\\displaystyle
\e\frac{\partial}{\partial t} v^\e(t,x)\! =\!(\frac{\partial^2}{\partial x^2}- m) v^\e(t,x)\! + \!
f(u^\e(t,x),v^\e(t,x))+  \rho(x) r(\alpha(t,x)) + 
 \e ^{1/2}\rho(x) \frac{\partial}{\partial t}
{\mathcal{W}}^2(t,x),   \\ \\
u^\e(t,0)=u^\e(t,1)=v^\e(t,0)=v^\e(t,1)=0, \,\\ \\
u^\e(0,x)= u^0(x), \ v^\e(0,x)= v^0(x),\qquad\qquad\qquad\qquad  t \in [0,1], \ x \in [0,1],
\end{array}\right.
\end{equation}
where $({\mathcal{W}}^1(t,x))$ and
 $({\mathcal{W}}^2(t,x))$  are independent
space-time white noises.  Here $(u^{\epsilon})$ represents the slow state, $(v^{\epsilon})$ the quick one and $\alpha$ is the control.

We make the following assumptions on the coefficients:
\begin{enumerate}
\item  $m$ is a positive constant.
\item  $b, f$ are  continuous maps, $b$ is bounded and Lipschitz continuous w.r.t to the first two variables uniformly w.r.t. the control, moreover
$f$  is  Lipschitz  continuous with a constant  smaller then $m$. 
\item $\sigma, \rho$ are measurable and bounded functions $ [0,1] \to \R$. Moreover we ask that $|\sigma (x)| \geq c_\sigma $, for a.e. $x \in [0,1]$ and a suitable constant $c_\sigma>0$.
\item  $r:  \R \to \R$ is a measurable and bounded map.
\item An admissibile control $\alpha$ is any bounded predictable  process $\alpha: \Omega \times [0,1] \times [0,1] \to \R$ and the cost functional is
\[ J^\e(u_0,v_0)= \E \int_0^1   \int_0^1  \ell (u^\e (t,x),v^\e (t,x), \alpha(t,x) ) \, dx \, dt + \int_0^1 h(u^\e (1,x)) \, dx,\]
with $\ell$ and $h$ Lipschitz continuous and bounded functions.
\end{enumerate}

The abstract formulation in $H=K=L^2(0,1)$ and $U=  L^2(0,1)$ is identical to the one in  \cite[section 5]{FuTes_BE}. In this same place  it is shown shown that Hypotheses  \ref{A.1}---\ref{A.5}, \ref{B.3}, \ref{B.4} and \ref{C.1} hold. Notice that  thus Theorem \ref{main} and Theorem \ref{main_controllo} apply.} 
%Let us set $H=K= U= L^2([0,1])$. Then
%$\frac{\partial}{ \partial x}{\mathcal{W}}^u(t,\cdot) = W^u(t)$,
%$\frac{\partial}{ \partial x}{\mathcal{W}}^v(t,\cdot) = W^v(t)$ are cylindrical Wiener processes in
%$ \Xi= L^2([0,1])$, see for instance \cite{DPZ1}. \\
%Moreover, we introduce $D(A)=D(B)= H^1_0([0,1]) \cap H^2([0,1])$ and the operators
%\[( A \xi) (x) =  \frac{\partial^2}{\partial x^2} \xi(x), \quad \xi \in D(A),\  \text{ and } \
%( B \xi) (x) =  (\frac{\partial^2}{\partial x^2}- m) \xi(x), \quad \xi \in D(B).  \]
%Finally, we put $ b(\xi,\zeta,\delta)(x)= b(\xi(x),\zeta(x), \delta(x)), \ F(\xi,\zeta)(x)= f(\xi(x),\zeta(x)),$
%$ R \xi(x)= \sigma(x)\xi(x),$ $ G\xi(x)= \xi(x) \rho(x). $
\end{example}

\section{Control interpretation of the limit forward-backward system}
Since we were able to interpret the limit value function as the solution of a \textit{reduced} forward we can now hope to see it as the value function of a correspondingly 
 \textit{reduced} control problem.

\medskip
 
Most of our analysis in this section is based on the fact that $\lambda$ is concave with 
respect to $z$. In particular, by Fenchel-Moreau theorem (translated in the obvious way for concave functions instead than for convex ones),
we can write $\lambda=\lambda_{**}$ where for all $x\in H$:
\begin{equation*}
\lambda_*(x,{p}) = \inf_{z \in \Xi^*} \bigl( -z p -\lambda(x,z)  \bigr) , \quad p\in \Xi
\end{equation*}
and the map $\lambda_*(x,\cdot)$ is an upper semicontinuous concave function with non empty domain in $\Xi$.
Thus for all $x\in H$, $z\in \Xi^*$:
\begin{equation*}
\lambda(x,z)=  \inf_{{p} \in \Xi} \bigl(- z {p} - \lambda_*(x,p)\bigr)
.\end{equation*}
Recalling that $\lambda$ is Lipschitz continuous with respect to $z$ uniformly in $x$ and denoting by 
$L$ the Lipschitz constant we have:
\begin{equation*}
\lambda_*(x,p) = -\infty,\hbox{ whenever } |p|>L.
\end{equation*}
and consequently:
\begin{equation*}
\lambda(x,z)=  \inf_{{p} \in \Xi, \, |p|\leq L} \bigl(- z {p} - \lambda_*(x,p)\bigr)
.\end{equation*}
 Moreover, $\lambda_*(x,{p}) \leq   -\lambda(x,0)\leq c(1+|x|)$, 
thus, for any process $(\mathfrak{p}_{t})_{0 \leq t \leq 1}$ with values in $\Xi$, the process 
\begin{equation*}
\biggl( \int_{0}^t \lambda_*(X_{s},\mathfrak{p}_{s}) ds \biggr)_{0 \leq t \leq 1}
\end{equation*}
is well-defined and takes values in $[-\infty,\infty)$. 

Given any $\Xi$ valued progressively measurable process $(\mathfrak{p}_t)_{t\geq 0}$ with $|\mathfrak{p}_t|\leq L$:
$$\begin{array}{rcl}
\bar{Y}_t &=& \displaystyle h(X_1)+\int_t^1 \lambda(X_s,\bar{Z}_s)ds -\int_t^1 \bar{Z}_sdW_s \\ &\leq& h(X_1) \displaystyle -\int_t^1 \left(\bar{Z}_s \mathfrak{p}_s+\lambda_*(X_s, \mathfrak{p}_s)\right) ds -\int_t^1 \bar{Z}_sdW_s.\end{array} $$
Introducing $W^{\mathfrak{p}}_t=\displaystyle \int_0^t \mathfrak{p}_s ds + W_t$ and the probability $\mathbb{P}^{\mathfrak{p}}$ under which it is a Wiener process we get:
$$dX_t=AX_t dt -R\mathfrak{p}_t dt +R dW^{\mathfrak{p}}_t,$$
$$\bar{Y}_t\leq h(X_1) -\int_t^1 \lambda_*(X_s, \mathfrak{p}_s) ds -\int_t^1 \bar{Z}_sdW^{\mathfrak{p}}_s,$$
%For such a process $(P_{t})_{0 \leq t \leq T}$, we deduce that 
%\begin{equation*}
%%&dX_{t} = AX_td{t} + dW^1_{t}
%%&d\bar{Y}_{t} \geq \bigl( Z_{t} P_{t} + \Lambda^*(X_{t},P_{t})\bigr) dt +  Z_{t} RdW^1_{t}, \quad Y_{1}=h(X_{1}),
%%\end{equation*}
which shows that: %, whenever $(\mathfrak{p}_{t})_{0 \leq t \leq 1}$ is bounded by $L$ 
%\begin{equation*}
%Y_{t} \leq \E_{t}^{\mathfrak{v}} \biggl[ %{\mathcal E}\biggl( \int_{0}^1R^{-1}P_{t} dW_{t}
%\biggr) \biggl( h(X_{1}) - \int_{0}^1 %\Lambda^*(X_{t},P_{t}) dt \biggr) \biggr],
%\end{equation*}
\begin{equation*}
\bar{Y}_{t} \leq \E^{\mathfrak{p}} \biggl( h(X_{1}) - \int_{t}^1 \lambda_*(X_{s},\mathfrak{p}_{s}) ds \bigg | \mathcal{F}_t\biggr).
\end{equation*}
%where we observe that the right-hand  side is always well-defined as an element 
%of $(-\infty,\infty]$.
% Above, we used the notation ${\mathcal E}$ for the Dol\'eans-Dade transformation. 
\vskip 4pt

Conversely, we may call, for any $n \geq 1$, 
$(\mathfrak{p}_{t}^n)_{0 \leq t \leq 1}$ such that %[CAMBIATO I SEGNI CONTROLLA!]
%$\lambda_*^\sharp(X_{t},\mathfrak{v}_{t}^n) \geq - \bar Z_{t} \mathfrak{v}_{t}^n -\lambda(X_{t},Z^n_{t}) + 1/n$ or equivalently
 $ - \bar Z_{t}\mathfrak{p}_{t}^n - \lambda_*(X_{t},\mathfrak{p}_{t}^n)-1/n\leq \lambda (X_{t},\bar{Z}_{t})$. Clearly we have $|\mathfrak{p}_{t}^n | \leq L$. Using a measurable selection Theorem, see for instance Theorem 6.9.13 in \cite{Bog}, one can choose the process $\mathfrak{p}^n$ to be progressive measurable. 

Then, we have 
$$\bar{Y}_t\geq  h(X_1)  -\int_t^1 \left(\bar{Z}_s \mathfrak{p}^n_s+\lambda_*(X_s, \mathfrak{p}^n_s)+\frac{1}{n}\right) ds -\int_t^1 \bar{Z}_sdW_s,
$$ 
and rewriting the above in terms of $W^{\mathfrak{p}^n}$:
\begin{equation*}
\begin{split}
\bar{Y}_{t} +\frac{1-t}{n} &\geq h(X_1)- \int_t^1  \lambda_*(X_{t},\mathfrak{p}_{t}^n) dt  - \int_t^1\bar Z_{t} dW^{\mathfrak{p}^n}_{t}.
\end{split}
\end{equation*}
Therefore we can conclude that $\bar{Y}_{t}$ is the value function of a stochastic optimal control problem in the sense that:
\begin{equation*}
\bar{Y}_{t} = \inf_{\mathfrak{p}} \E^{\mathfrak{p}}  \biggl( h(X_{1}) - \int_{0}^1 \lambda_*(X_{t},\mathfrak{p}_{t}) dt \bigg | \mathcal{F}_t \biggr),
\end{equation*}
where $(X_t)_{t\geq 0} $ is the solution of the following controlled
stochastic differential equation: 
$$dX_t=AX_t dt -R\mathfrak{p}_t dt +R dW^{\mathfrak{v}}_t, \; X_0=x_0,$$
 the supremum is extended to all $\Xi$-valued, predictable processes 
$(\mathfrak{p}_s)_{0 \leq s \leq 1}$ that are bounded by $L$ and finally $W^{\mathfrak{v}}$ is a $\Xi$-valued Wiener process with respect to a certain $\mathbb{P}^{\mathfrak{v}}$.

\medskip

\textbf{Acknowledgements:}  The authors wish to thank Francois Delarue for several 
enlightening discussions on the content of the present paper.

%%%%%%%%%%%%%%%%%%%%%%%%%%%%%%%%%%%%%%%%%%%%%%%%%%%%%%%

\end{document}